\documentclass{article}
\usepackage{amsfonts}
\usepackage{amsmath}
\usepackage{amssymb}
\usepackage{algorithm}
\usepackage{algorithmic}
\usepackage[english]{babel}
\usepackage{booktabs}
\usepackage{caption}
\usepackage{colortbl}
\usepackage{etoolbox}
\usepackage[T1]{fontenc}
\usepackage{float}
\usepackage{fancyhdr}
\usepackage[a4paper, portrait, margin=1.1811in]{geometry}
\usepackage{graphicx}
\usepackage{helvet}
\usepackage[colorlinks, citecolor=cyan]{hyperref}
\usepackage[utf8]{inputenc}
\usepackage{mathrsfs}
\usepackage{multirow}
\usepackage[thmmarks,amsmath]{ntheorem}
\usepackage{pgfplots}
\usepackage{scrextend}
\usepackage{subfigure}%
\usepackage{titlesec}
\usepackage{xcolor}

\pgfplotsset{compat=1.15}

\usetikzlibrary{arrows}
\graphicspath{ {./images/} }
\newcounter{lemma}
\newtheorem{lemma}{Lemma}
\newcounter{theorem}
\newtheorem{theorem}{Theorem}
\newtheorem{remark}{Remark}
\newcounter{definition}
\newtheorem{definition}{Definition}
\newtheorem{proposition}{Proposition}

\newcounter{example}
\theorembodyfont{\normalfont}
\newtheorem{example}{Example}
\theoremstyle{nonumberplain}
\theoremheaderfont{\bfseries}
\theorembodyfont{\normalfont}
\theoremsymbol{$\square$}
\newtheorem{Proof}{\hskip 0em Proof}

\makeatletter
\patchcmd{\@maketitle}{\LARGE \@title}{\fontsize{16}{19.2}\selectfont\@title}{}{}
\makeatother

\usepackage{authblk}

\newsavebox\affbox

\author[1]{Shicong Zhong}
\author[1]{Bingru Huang\footnote{Corresponding author. E-mail: \href{}{hbr999@ustc.edu.cn}}}
\author[1]{Falai Chen}
\affil[1]{School of Mathematical Sciences, University of Science and Technology of China, Hefei, Anhui 230026, P. R. China}

\titlespacing\section{0pt}{12pt plus 4pt minus 2pt}{0pt plus 2pt minus 2pt}
\titlespacing\subsection{12pt}{12pt plus 4pt minus 2pt}{0pt plus 2pt minus 2pt}
\titlespacing\subsubsection{12pt}{12pt plus 4pt minus 2pt}{0pt plus 2pt minus 2pt}

\titleformat{\section}{\normalfont\fontsize{10}{15}\bfseries}{\thesection.}{1em}{}
\titleformat{\subsection}{\normalfont\fontsize{10}{15}\bfseries}{\thesubsection.}{1em}{}
\titleformat{\subsubsection}{\normalfont\fontsize{10}{15}\bfseries}{\thesubsubsection.}{1em}{}

\titleformat{\author}{\normalfont\fontsize{10}{15}\bfseries}{\thesection}{1em}{}

\title{\textbf{\huge Basis construction for polynomial spline spaces over arbitrary T-meshes}} 
\date{}    

\begin{document}

\pagestyle{headings}	
\newpage
\setcounter{page}{1}
\renewcommand{\thepage}{\arabic{page}}

\captionsetup[figure]{labelfont={bf},labelformat={default},labelsep=period,name={Figure }}	\captionsetup[table]{labelfont={bf},labelformat={default},labelsep=period,name={Table }}
\setlength{\parskip}{0.5em}
	
\maketitle
	
\noindent\rule{15cm}{0.5pt}
	\begin{abstract}
		This paper presents a novel method for constructing bases for  polynomial spline spaces  over \textcolor{blue}{arbitrary T-meshes} (PT-splines for short). We construct spline basis functions for an arbitrary T-mesh  by first converting the T-mesh into a diagonalizable one via edge extension, ensuring a stable dimension of the spline space. Basis functions over the diagoalizable T-mesh are constructed 
        according to the three components in the dimension formula corresponding to cross-cuts, rays, and T $l$-edges in the diagonalizable T-mesh, and each component is assigned some local tensor product B-splines as the basis functions. We prove this set of functions constitutes a basis for the diagonalizable T-mesh. To remove redundant edges from extension, we introduce a technique, termed Extended Edge Elimination (EEE) to construct a basis for an arbitrary T-mesh while reducing structural constraints and unnecessary refinements. The resulting PT-spline basis ensures linear independence and completeness, supported by a dedicated construction algorithm. A comparison with LR B-splines, which may lack linear independence and are limited to LR-meshes, highlights the PT-spline's versatility across any T-mesh. Examples are also provided to demonstrate that dimensional instability in spline spaces
        is related with 
        basis function degradation and that PT-splines are advantageous over HB-splines for certain hierarchical T-meshes.\\  \\
		\textbf{\textit{Keywords}}: \textit{Polynomial splines over arbitrary T-meshes;
        PT-spline basis; Dimensional instability; Isogeometric analysis}
	\end{abstract}
\noindent\rule{15cm}{0.4pt}

\section{Introduction}
Splines have found widespread applications across numerous fields, including computer-aided geometric design~\cite{lin1996nurbs,tolvaly2015mixed,mykhaskiv2018nurbs}, engineering~\cite{ma1998nurbs,yin2004reverse}, scientific computing~\cite{abdellah2025high}, the robotics \cite{yan2015representation,wu2016path}, and even biomedical sciences\cite{verhaeghe2007reconstruction,chen2010parametric}. The study of spline functions began in the mid-20th century, originally as an alternative to polynomials for interpolation and function approximation~\cite{schoenberg1946contributions}. In the 1960s, with the rise of computers, splines became an integral part of computer-aided design, revolutionizing shape modeling and industrial manufacturing~\cite{ferguson1964multivariable}. De Boor's extensive research on B-splines ~\cite{de1978practical} played a key role in refining their theoretical framework and establishing them as a fundamental mathematical tool. 
A major breakthrough came in 1975 with the introduction of Non-Uniform Rational B-Splines (NURBS), which became the standard for multivariate spline modeling. 

Since NURBS are defined on tensor-product grids, their mathematical properties are relatively straightforward to derive. However, this also limits their flexibility, as they do not support local refinement, often resulting in an excessive number of redundant control points. To overcome this limitation, alternative spline constructions on T-meshes, such as T-splines~\cite{sederberg2003t}, LR-splines~\cite{dokken2013polynomial}, PHT-splines~\cite{deng2008polynomial} and (T)HB-splines~\cite{forsey1988hierarchical, giannelli2012thb}, have been developed. T-splines introduce T-junctions in the control mesh, allowing control points to be inserted without affecting entire rows or columns and thus tackle the limitations of NURBS~\cite{sederberg2004t,sederberg2008watertight}. T-splines is applied widely in geometric modeling and isogeometric analysis (IGA)~\cite{bazilevs2006isogeometric,bazilevs2010isogeometric,dorfel2010adaptive,scott2011isogeometric}. However, they do not always preserve linear independence and analysis-suitable T-splines (AST splines) is proposed in ~\cite{buffa2010linear,li2012linear,scott2012local,li2014analysis} to overcome the limitation. LR-splines, are constructed through the refinement of spline functions and form a set of functions within the polynomial spline space. While they offer adaptive refinement, their subdivision is restricted by the mesh structure, and they do not always guarantee linear independence, a crucial property in the numerical analysis of differential equations. PHT splines are spline spaces with reduced regularity (such as $S(3,3,1,1)$)  which have perfect local refinement property and have been widely used in IGA~\cite{li2007surface,deng2008polynomial,nguyen2011isogeometric,wang2011adaptive}. HB-splines and THB-splines generate spline functions through hierarchical refinement, with the latter employing a truncation mechanism to enforce the partition of unity property. However, as noted in ~\cite{mokrivs2014completeness}, although they maintain linear independence, additional constraints on the mesh topology are required to ensure they serve as valid basis functions. 

The spline space over T-mesh was first introduced in \cite{deng2006dimensions}, and its study can be broadly divided into two key aspects: the exploration of basis functions and the analysis of dimensionality. Various methods have been developed to compute the dimension of spline spaces, including the B-net method~\cite{deng2006dimensions}, the smoothing cofactor-conformality method~\cite{wang1975structural,renhong2001multivariate}, and homological techniques~\cite{mourrain2014dimension}.
In \cite{deng2006dimensions}, the B-net method was employed to derive a dimension formula for the spline space $S(d_1,d_2,\alpha,\beta)$,with $d_1\geq 2\alpha+1; d_2\geq 2\beta+1$. The smoothing cofactor-conformality method offers a way to express smoothness conditions in algebraic form. In ~\cite{huang2006new}, this method was applied to provide an alternative proof of the dimension formula for the same space as in ~\cite{deng2006dimensions}. Additionally, ~\cite{li2011instability} used the same approach to reveal the instability of spline spaces under conditions of maximal smoothness. Moreover, ~\cite{guo2015problem} conducted an in-depth investigation into the instability of spline space dimensions when the mesh features a specific topological structure known as a T-cycle. However, despite their focus on dimensionality, these studies generally do not establish a direct connection between dimension analysis and the construction of basis functions except PHT-splines.

This study presents a method for constructing spline basis functions over \textbf{arbitrary T-meshes}. We  first convert an arbitrary T-mesh into a diagonalizable T-mesh through edge extension, the spline space over which has a stable dimension.  The dimension formula of the spline space consists of three parts which correspond to the three components of the diagonalizable T-mesh: cross-cuts, rays, and T $l$-edges. We then assign some local tensor product B-splines to each component. We prove that this set of functions forms a basis for the spline space over the diagonalizable T-mesh. In order to remove the auxiliary edged introduced by edge extension, we propose a method, termed \textbf{Extended Edge Elimination} (EEE for short) to produce a basis for the spline space over the given T-mesh. This approach reduces mesh structure constraints and minimizes unnecessary refinements during basis construction.

The main contributions of this paper are as follows:
\begin{itemize}
    \item We develop a  method for  constructing a basis of the spline space over a diagonalizable T-mesh.
    \item We introduce the  Extended Edge Elimination, a method essential for constructing basis functions over an arbitrary T-mesh.
    \item We develop a method to construct basis functions over an arbitrary T-mesh and prove that these basis functions exhibit local supports.
\end{itemize}

This article is structured as follows. Section 2 provides a review of polynomial spline spaces over T-meshes, including some basic notions for T-meshes, diagonalizable T-meshes and \(t\)-partitions. Section 3 describes basis construction method for spline spaces over arbitrary T-meshes. In Section 4, we present more examples of the PT-spline basis construction, compare it with LR B-splines, and offer a preliminary analysis of dimensional instability in spline spaces over T-meshes. Finally, Section 5 concludes the paper and proposes topics for future research.

\section{Preliminaries}
In this section, we introduce some preliminaries include T-mesh and spline space over it, smoothing cofactor-conformality method, diagonalizable T-mesh and the dimension formula of it.

\subsection{Spline spaces over arbitrary T-mesh}
First, we review the definition of T-mesh and the spline spaces over it.

\begin{definition}(\cite{deng2006dimensions, huang2024stability}) Suppose $\mathscr{T}:=\{C_i\}_{i\in I}$ be a collection of axis-aligned rectangle boxes and the domain occupied by these boxes is $\Omega(\mathscr{T}):=\cup_{i\in I} C_i$ where $I$ is the index set. In addition, if any pair of distinct rectangles can only intersect at points on their edges and $\Omega(\mathscr{T})$ is connected, $\mathscr{T}$ is called a \textbf{T-mesh} and $C_i, i\in I$ is called a \textbf{cell} of $\mathscr{T}$. We call \textcolor{blue}{ a T-mesh} $\mathscr{T}$ a \textbf{regular T-mesh} if $\Omega (\mathscr{T})$ has no holes and its boundary is a rectangle. 
\end{definition}

We give some supplementary explanations for certain terms used in the T-mesh. A grid point in a T-mesh is called a \textbf{vertex} of the T-mesh. If a vertex is on the boundary of T-mesh domain, it is called a \textbf{boundary vertex}, otherwise it is called a \textbf{interior vertex}. Note that an interior vertex can be a \textbf{T-node} which is the intersection of three cells in a T-mesh. An interior vertex that is the intersection of four cells is a \textbf{crossing node}. The line segment connecting two \emph{adjacent} vertices on a grid line together with the two vertices is called an \textbf{edge} of a T-mesh. A \textbf{$l$-edge} is the longest possible line segment containing several edges and its end points are T-nodes or boundary vertices. There are three types of $l$-edges in a T-mesh: 1. A \textbf{T $l$-edge} is the $l$-edge whose two end points are both T-nodes; 2. \textcolor{blue}{A \textbf{ray line} }is the $l$-edge whose one end point is T-node and another is boundary point; 3. \textcolor{blue}{A \textbf{cross-cut}} is the $l$-edge whose two end points are both boundary points. We introduce some more notations in Table \ref{notation} for convenience and use two end points to note a segment line in the subsequent discussion. 

\begin{figure}[htbp]  
    \centering      %
    \includegraphics[width=0.4\textwidth]{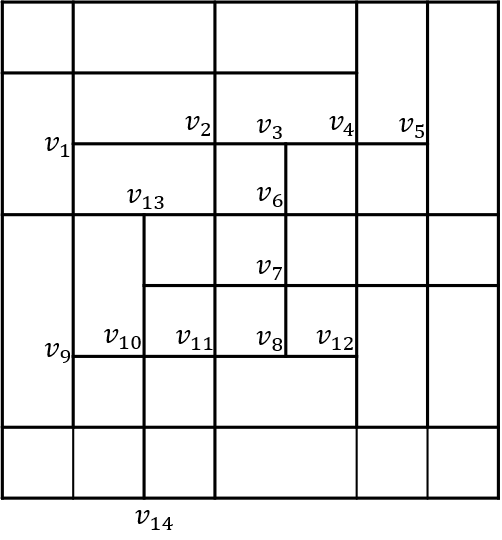} 
    \captionsetup{font={scriptsize}}
    \caption{T-mesh} 
    \label{wholeex4}  %
\end{figure}

\begin{example}
    Figure \ref{wholeex4} shows a regular T-mesh, \textcolor{blue}{$v_1,v_3,v_5,v_8,v_9,v_{12},v_{13}$ are T-nodes, $v_2,v_4,v_6,v_7,v_{10},v_{11}$ are crossing nodes} and $v_{14}$ is boundary vertex. $v_{13}v_{14}$ is a ray line and $v_1v_5,v_3v_8,v_9v_{12}$ are T $l$-edges. 
\end{example}

\begin{table}[htbp]
\centering
\caption{Notations for T-mesh $\mathscr{T}$}
\setlength{\tabcolsep}{15pt}
\renewcommand{\arraystretch}{2}
\begin{tabular}{cc}
\toprule
$c_h$& Number of horizontal cross-cuts in T-mesh\\
$c_v$& Number of vertical cross-cuts in T-mesh\\
$t_h$& Number of horizontal T $l$-edges in T-mesh\\
$t_v$& Number of vertical T $l$-edges in T-mesh\\
$t$& Number of all T $l$-edges in T-mesh\\
$n_v$  & Number of interior vertices in T-mesh\\
$n_T$& Number of interior vertices on T $l$-edges T-mesh\\
$T(\mathscr{T})$ & The set of all T $l$-edges \\
\bottomrule
\end{tabular}
\label{notation}
\end{table}

\begin{definition}(\cite{deng2006dimensions}) Let $\mathscr{T}$ be a T-mesh, with $\Omega$ denoting the region covered by its cells. The spline space over $\mathscr{T}$ is defined as
\[
S(d_1,d_2,\alpha,\beta,\mathscr{T}) := \{ s(x,y) \in C^{\alpha,\beta}(\Omega) \mid s(x,y)|_{C_i} \in \mathbb{P}_{d_1,d_2}, \ \forall C_i \in \mathscr{T} \},
\]
where $\mathbb{P}_{d_1,d_2}$ is the space of bivariate polynomials of degree $(d_1,d_2)$, and $C^{\alpha,\beta}(\Omega)$ is the space of bivariate functions continuous in $\Omega$ with order $\alpha$ in the $x$-direction and order $\beta$ in the $y$-direction. The space $S(d_1,d_2,\alpha,\beta,\mathscr{T})$ is linear. For the highest smoothness case, we denote $S_{d_1,d_2}(\mathscr{T}) := S(d_1, d_2, d_1-1, d_2-1, \mathscr{T})$.
\end{definition}

An T $l$-edge $l$ is defined as a \textbf{vanished $l$-edge} if it contains at most $d_1+1$ vertices when horizontal or at most $d_2+1$ vertices when vertical.

\subsection{Diagonalizable T-mesh and $t$-partition}

This section briefly reviews the concepts of diagonalizable T-mesh~\cite{li2016dimension} and \( t \)-partition~\cite{huang2025}, which are utilized throughout the paper.

For a given T-mesh \(\mathscr{T}\) with ordered T \(l\)-edges \(l_1 \succ l_2 \succ \cdots \succ l_{t}\), we denote by \(\overline{l}_i\) the T \(l\)-edge \(l_i\) with all intersections with \(l_j\) (\(j=1, 2, \ldots, i-1\)) removed.

\begin{definition}(\cite{huang2025})
    Given a T-mesh \(\mathscr{T}\) with T \(l\)-edges ordered as \(l_1 \succ l_2 \succ \cdots \succ l_t\), the set \(\{\overline{l}_1, \overline{l}_2, \ldots, \overline{l}_t\}\) forms a partition of \(T(\mathscr{T})\), termed a \textbf{\(t\)-partition}. Any edge \(\overline{l}_i,i=1,2,\ldots,t\) is defined as \textbf{an \(l\)-edge in a \(t\)-partition}, and we consider \(\overline{l}_i,i=1,2,\ldots,t\) as a specific type of \(l\)-edge.
\end{definition}

Note that the $t$-partition is a concept independent of any ordering. And $$\sum\limits_{i=1}^{t}n(\overline{l}_i)=n_{T}.$$ where $n(\overline{l}_i)$ is the number of vertices on $\overline{l}_i$.

\begin{figure}[htbp]  
    \centering      %
    \includegraphics[width=0.4\textwidth]{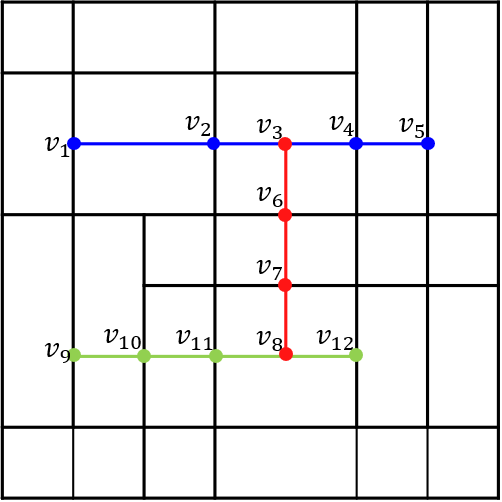} 
    \captionsetup{font={scriptsize}}
    \caption{t-partition} 
    \label{wholeex1}  %
\end{figure}

\begin{example}\label{tpart}
    In Figure~\ref{wholeex1}, there are three T \(l\)-edges: \(l_1 = v_3v_8\) marked in red, \(l_2 = v_9v_{12}\) marked in light green, and \(l_3 = v_1v_5\) marked in blue. These can be ordered as \(l_1 \succ l_2 \succ l_3\), where \(\overline{l}_1 = l_1\), \(\overline{l}_2 = l_2\) with \(v_8\) removed, and \(\overline{l}_3 = l_3\) with \(v_3\) removed. The set \(\{\overline{l}_1, \overline{l}_2, \overline{l}_3\}\) constitutes a $3$-partition of this T-mesh.
\end{example}

\begin{definition}(\cite{li2016dimension,huang2025})
A T-mesh \(\mathscr{T}\) is \textbf{diagonalizable} if there exists a \(t\)-partition \(\{\overline{l}_1, \overline{l}_2, \ldots, \overline{l}_{t}\}\) such that \(n(\overline{l}_i) \geq d_1 + 1, i=1,2,\ldots,t_h\) for horizontal \(l_i\) or \(n(\overline{l}_j) \geq d_2 + 1,j=1,2,\ldots,t_v\) for vertical \(l_j\), where \(n(\overline{l}_i)\) denotes the number of vertices on \(\overline{l}_i\).
\end{definition}

In \cite{li2016dimension}, Li et al. derive an explicit dimension formula for the spline space over a T-mesh:

\begin{theorem}(\cite{li2016dimension})
For a diagonalizable T-mesh $\mathscr{T}$ which has no vanished $l$-edges and holes, it's dimension of the spline space $S_{d_1,d_2}(\mathscr{T})$ is 
\begin{equation}\label{dim diag}
    \dim S_{d_1,d_2}(\mathscr{T})=(d_1+1)(d_2+1)+c_h(d_1+1)+c_v(d_2+1)+n_v-((d_1+1)t_h+(d_2+1)t_v).
\end{equation}
\end{theorem}

In \cite{huang2025}, Huang et al. prove that, in a diagonalizable T-mesh, each T $l$-edge in the $t$-partition corresponds to a one-dimensional B-spline, and these edges are mutually independent. In the subsequent sections, we demonstrate that both T $l$-edges and rays in a diagonalizable T-mesh correspond to one-dimensional B-splines, while cross-cuts correspond to two-dimensional tensor product B-splines.

\section{The construction of spline basis over an arbitrary T-mesh}
\textcolor{blue}{In this section, we introduce a novel method for constructing basis functions over an arbitrary T-mesh. For any given T-mesh, we first modify it to a diagonalizable T-mesh by edge extension. Then we construct a set of basis functions which are local tensor-product B-splines on the extended diagonalizable T-mesh, which is similar to T-spline construction~\cite{sederberg2003t}. Finally, we  propose an Extended Edge Elimination technique to construct a set of basis functions over the original T-mesh. The resulting basis functions consist of local tensor-product B-splines and their linear combinations, similar to those in THB-splines~\cite{giannelli2012thb}.
Several key properties of the basis functions are also established.}


\subsection{Basis of spline space over extended T-Meshes}\label{subsection 3.1}

\begin{definition}
  Let $\mathscr{T}$ be a given T-mesh, we define \textbf{the multiplicity function of an edge} as 
  \[
\mu(l) =
\begin{cases}
1 & \text{if } \text{$l$ is an interior edges of $\mathscr{T}$}, \\
d_1+1 & \text{if } \text{$l$ is one of the vertical boundary edges of $\mathscr{T}$}\\
d_2+1 & \text{if } \text{$l$ is one of the horizontal boundary edges of $\mathscr{T}$}.
\end{cases}
\]
\end{definition}

In a specified coordinate system, let $\mathsf{X}$ and $\mathsf{Y}$ represent the sets of all $x$-coordinates and $y$-coordinates, respectively, of the knots within the T-mesh. Subsequently, we define the local tensor product B-spline associated with the $l$-edge of $\mathscr{T}$ as follows.

\begin{definition}\label{def local B associated}
    Given a T-mesh \(\mathscr{T}\), for a horizontal \(l\)-edge \(l\), let \(v_1, v_2, \ldots, v_{d_1+2}\) be any \(d_1+2\) vertices on \(l\) with corresponding \(x\)-coordinates \(x_{\sigma_1}, x_{\sigma_2}, \ldots, x_{\sigma_{d_1+2}} \in \mathsf{X}\). We define the local tensor product B-spline associated with \(l\) at vertices \(v_1, v_2, \ldots, v_{d_1+2}\) as
\[
B[v_1, v_2, \ldots, v_{d_1+2}](x, y) := N(x_{\sigma_1}, x_{\sigma_2}, \ldots, x_{\sigma_{d_1+2}})(x) \times N(y_{\tau_1}, y_{\tau_2}, \ldots, y_{\tau_{d_2+2}})(y),
\]
where \(N(x_{\sigma_1}, x_{\sigma_2}, \ldots, x_{\sigma_{d_1+2}})(x)\) is the one-dimensional B-spline defined by the nodes \(\{x_{\sigma_i}\}_{i=1}^{d_1+2}\), provided there exist \(y_{\tau_1}, y_{\tau_2}, \ldots, y_{\tau_{d_2+2}} \in \mathsf{Y}\) such that the line segments \(\{[x_{\sigma_1}, x_{\sigma_{d_1+2}}] \times y_{\tau_i}\}_{i=1}^{d_2+2}\) and \(\{x_{\sigma_i} \times [y_{\tau_1}, y_{\tau_{d_2+2}}]\}_{i=1}^{d_1+2}\) lie within the \(l\)-edges of \(\mathscr{T}\), with one segment coinciding with \(l\). Here, \(N(y_{\tau_1}, y_{\tau_2}, \ldots, y_{\tau_{d_2+2}})(y)\) denotes the one-dimensional B-spline defined by the nodes \(\{y_{\tau_j}\}_{j=1}^{d_2+2}\). If these conditions are not satisfied, the local tensor product B-splines or local tensor product structure are considered not associated with $l$ at vertices \(v_1, v_2, \ldots, v_{d_1+2}\). 

The line segments \(\{x_{\sigma_i} \times [y_{\tau_1}, y_{\tau_{d_2+2}}], [x_{\sigma_1}, x_{\sigma_{d_1+2}}] \times y_{\tau_j}\}_{i=1,j=1}^{d_1+2,d_2+2}\) are defined as \textbf{the local tensor product structure associated with \(l\) at vertices \(v_1, v_2, \ldots, v_{d_1+2}\)}. 

These notions for a vertical \(l\)-edge can be defined analogously.

\end{definition}

\begin{figure}[htbp]  
    \centering      %
    \includegraphics[width=0.4\textwidth]{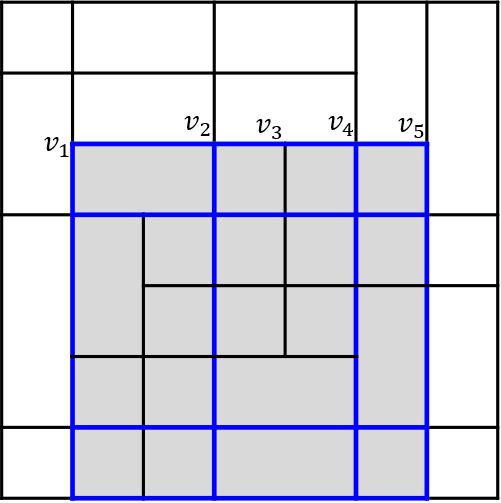} 
    \captionsetup{font={scriptsize}}
    \caption{The tensor product structure associated with $v_1v_5$ at vertices $v_1,v_2,v_4,v_5$.} 
    \label{wholeex2}  %
\end{figure}

\begin{example}
    Consider the spline space \(S_{2,2}(\mathscr{T})\); the gray domain in Figure~\ref{wholeex2} represents the support of a tensor product structure, which is one of the tensor product structures associated with the T \(l\)-edge \(v_4v_5\) at vertices \(v_1, v_2, v_4, v_5\).
\end{example}

Consider a diagonalizable T-mesh. Utilizing the dimension formula \eqref{dim diag}, we derive a simplified expression as follows:
\begin{equation}\label{dim diag simplify}
    \left((d_1+1+c_h)(d_2+1+c_v)\right)+ \left(n_T - ((d_1+1)t_h+(d_2+1)t_v)\right) + \left(n_v-n_T-c_vc_h\right).
\end{equation} 
The simplified dimension formula \eqref{dim diag simplify} guides the partition of the diagonalizable T-mesh into three components:
\begin{itemize}
\item The tensor product component (cross-cuts and boundary lines);
\item The T $l$-edges component within the $t$-partition;
\item The ray lines component.
\end{itemize}
These components align with the terms in the dimension formula:
\begin{itemize}
\item The dimension contributed by the tensor product component is
\begin{equation}\label{eq cross-cut}
    \textcolor{blue}{(d_1+1+c_h)(d_2+1+c_v).}
\end{equation}   
\item The dimension contributed by the T $l$-edges component is
\begin{equation}\label{eq T l-edges}
    n_T - ((d_1+1)t_h+(d_2+1)t_v).
\end{equation}  
\item The dimension contributed by the ray lines component is
\begin{equation}\label{eq ray}
    n_v-n_T-c_vc_h.
\end{equation}  

\end{itemize}
We define the T-mesh $\mathscr{T}$ with its rays and T $l$-edges removed as the \textbf{tensor product part of $\mathscr{T}$}. It is clear that, for each cross-cut in $\mathscr{T}$, a local tensor product B-spline associated with the cross-cut can be identified within the tensor product part of $\mathscr{T}$ using consecutive nodes in both $x$ and $y$ directions, with the total number of such B-splines corresponding to the dimension given by \eqref{eq cross-cut}. Consequently, the subsequent analysis will focus solely on T $l$-edges and rays.

\begin{figure}[htbp]  
    \centering      %
    \includegraphics[width=0.4\textwidth]{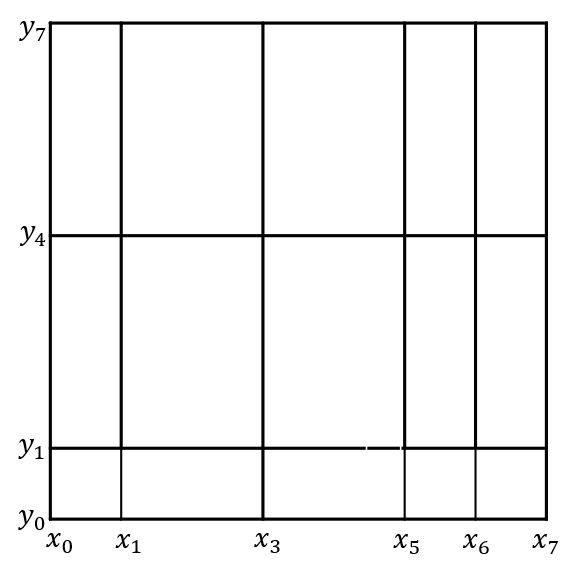} 
    \captionsetup{font={scriptsize}}
    \caption{The tensor product part of the T-mesh in Figure \ref{wholeex4}} 
    \label{wholeex6}  %
\end{figure}

\begin{example}
    Figure \ref{wholeex6} is the tensor product part of T-mesh shown in Figure \ref{wholeex4}, and for this part we define 35 local tensor product B-splines which are $$\{N(x_{\sigma_i},x_{\sigma_{i+1}},x_{\sigma_{i+2}},x_{\sigma_{i+3}})(x)\times N(y_{\sigma_j},y_{\sigma_{j+1}},y_{\sigma_{j+2}},y_{\sigma_{j+3}})(y)\}^{6,4}_{i,j=0}$$ where $(x_{\sigma_0},...,x_{\sigma_9})=(x_0,x_0,x_0,x_1,x_3,x_5,x_6,x_7,x_7,x_7)$, $(y_{\sigma_0},...,y_{\sigma_7})=(y_0,y_0,y_0,y_1,y_4,y_7,y_7,y_7)$ in the polynomial spline space $S_{2,2}(\mathscr{T})$.
\end{example}

Given that not every T $l$-edge or ray in the T-mesh $\mathscr{T}$ possesses the necessary associated local tensor product B-spline, as defined in Definition~\ref{def local B associated}, we proceed to investigate the extension of the T-mesh as follows.

\begin{definition}
    Given a T-mesh \(\mathscr{T}\), the T-mesh \(\mathscr{T}_1\) is defined as an extension of \(\mathscr{T}\) if \(\mathscr{T}_1\) is derived by extending the rays or T \(l\)-edges of \(\mathscr{T}\), denoted as \(\mathscr{T}\subseteq\mathscr{T}_1 = \mathrm{ext}_s(\mathscr{T})\), where \(s\) represents the number of edges extended from \(\mathscr{T}\) to \(\mathscr{T}_1\).
\end{definition}

Obviously we have the following ascending chain:
$$\mathrm{ext}_0(\mathscr{T})=\mathscr{T}\subseteq \mathrm{ext}_1(\mathscr{T})\subseteq \mathrm{ext}_2(\mathscr{T})\subseteq\cdots\subseteq\mathrm{ext}_{r}(\mathscr{T}).$$
where $\mathrm{ext}_{r}(\mathscr{T})$ is a tensor product mesh called \textbf{the extended tensor product mesh of $\mathscr{T}$.}

In the following, we demonstrate that for any T-mesh \(\mathscr{T}\), there exists an extended T-mesh \(\mathrm{ext}_s(\mathscr{T})\) such that \(\mathrm{ext}_s(\mathscr{T})\) is a diagonalizable T-mesh and ensures that the number of the local tensor product B-splines associated with T \(l\)-edges and rays correspond to the dimension formula of the respective components. Prior to presenting the main theorem, we introduce several preparatory lemmas.

\begin{lemma}\label{lem one dimension indp}
Suppose there are $N$ one-dimensional B-splines $B_i(x),i=1,2,\ldots,N$ of degree $d$, with the knots 
$$(\underbrace {x_0,...,x_0}_{d+1},x_i),$$
for $i=1,...,N$ where $x_0<x_1<x_2<\ldots<x_N$, then these $N$ one-dimensional B-splines are \textcolor{blue}{linearly independent }on the interval $[x_0,x_N]$.
\end{lemma}
\begin{Proof}
    Suppose \begin{equation}\label{eq one dimensional linear independent}
        \sum_{i=1}^N c_i B_i(x) = 0
    \end{equation} for \(x \in [x_0, x_N]\), where 
    $$B_i(x) = N(\underbrace{x_0, \ldots, x_0}_{d+1}, x_i)(x).$$ 
    We aim to prove that \(c_i = 0\) for \(i = 1, 2, \ldots, N\).

    Consider the equation~\eqref{eq one dimensional linear independent} restricted to the interval \((x_{N-1}, x_N)\), where only \(B_N(x)\) is nonzero. This reduces the equation to \(c_N B_N(x) = 0\). Since \(B_N(x) \neq 0\) for \(x \in (x_{N-1}, x_N)\), it follows that \(c_N = 0\). Extending this process to all intervals \((x_{i-1}, x_i)\) for \(i\) from \(N-1\) down to 1, we obtain \(c_i = 0\) for \(i = 1, \ldots, N\). Thus, the functions \(\{B_i(x)\}_{i=1}^N\) are linearly independent.
\end{Proof}

For the case where \(x_N < x_{N-1} < \cdots < x_1 < x_0\), similar conclusions can be derived using the same approach, and thus, detailed elaboration is omitted here.

\begin{lemma}\label{lem existence}
 Given a T-mesh \(\mathscr{T}\), we consider its extended T-mesh \(\mathrm{ext}_s(\mathscr{T})\). Assume there are \(t'\) T \(l\)-edges within \(\mathrm{ext}_s(\mathscr{T})\). For each \(l\)-edge \(\overline{l}_i\) (\(i = 1, 2, \ldots, t'\)) in a \(t'\)-partition, let \(n(\overline{l}_i)\) denote the number of vertices on \(\overline{l}_i\), and \textcolor{blue}{let $V_p$ denote a set of $d_k + 2$ vertices $v_1^{(p)}, v_2^{(p)}, \ldots, v_{d_k+2}^{(p)}$ lying on $\overline{l}_i$, associated with the $p$-th local tensor-product B-spline whose support overlaps with (a portion of) $\overline{l}_i$}, where \(k = 1\) for horizontal \(l\)-edges and \(k = 2\) for vertical \(l\)-edges. Then, there exists an integer \(0 \le s < r\) such that the extended T-mesh \(\mathrm{ext}_s(\mathscr{T})\) forms a diagonalizable T-mesh, satisfying the following condition:
\begin{itemize}
    \item It includes at least \(n(\overline{l}_i) - d_k - 1\) local tensor product B-splines associated with each \(\overline{l}_i\) across \(n(\overline{l}_i) - d_k - 1\) distinct vertex sets \(V_1, \ldots, V_{n(\overline{l}_i) - d_k - 1}\), and for all \(1 \le j < i\), the \(n(\overline{l}_i) - d_k - 1\) local tensor product B-splines associated with \(\overline{l}_i\) are not associated with \(\overline{l}_j\). Furthermore, it is possible to select suitable \(V_1, \ldots, V_{n(\overline{l}_i) - d_k - 1}\) such that the one-dimensional B-splines defined by the corresponding nodes are linearly independent.
\end{itemize}
\end{lemma}

\begin{Proof}
    Consider the extended T-mesh $\mathrm{ext}_{s}(\mathscr{T})$, formed by expanding all horizontal T $l$-edges and horizontal rays in $\mathscr{T}$ into cross-cuts. Since all vertical T $l$-edges are non-vanishing $l$-edges and remain non-intersecting, we obtain
$$n(\overline{l}_i) = n(l_i) \ge d_2 + 1, \quad i = 1, 2, \ldots, t_v,$$
thereby confirming that $\mathrm{ext}_{s}(\mathscr{T})$ constitutes a diagonalizable T-mesh.
Additionally, as the boundary of each horizontal cross-cut contains \(2(d_1 + 1)\) multiplicity knots, we select \(d_1 + 1\) boundary knots and one vertex in $\overline{l}_i$ in the \(x\)-direction for each \(\overline{l}_i\) and \(V_p\) for some index \(p\) in the \(y\)-direction. By Lemma~\ref{lem one dimension indp}, it is straightforward to verify that the construction of local tensor product B-splines associated with all T \(l\)-edges in the \(t_v\)-partition satisfies the conditions of the Lemma~\ref{lem existence}.
\end{Proof}

\begin{figure}[htbp]  
    \centering      %
    \includegraphics[width=0.4\textwidth]{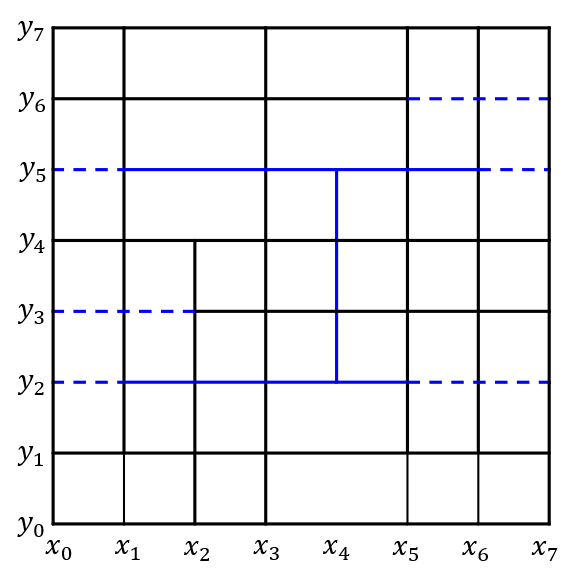} 
    \captionsetup{font={scriptsize}}
    \caption{Extended T-mesh as the proof of Lemma \ref{lem existence}} 
    \label{wholeex5}  %
\end{figure}

\begin{example}
  As illustrated in Figure \ref{wholeex5}, we extend all horizontal T $l$-edges and horizontal rays in Figure \ref{wholeex1} into cross-cuts. Then this extended T-mesh satisfies the condition in Lemma \ref{lem existence}.
\end{example}

Lemma~\ref{lem existence} demonstrates that an extended T-mesh $\mathrm{ext}(\mathscr{T})$ can be constructed by extending the T $l$-edges and rays of $\mathscr{T}$, satisfying the following condition:
\begin{itemize}
\item For all T $l$-edges in a $t'$-partition of $\mathrm{ext}_s(\mathscr{T})$ where $t'$ is the number of T $l$-edges in $\mathrm{ext}_s(\mathscr{T})$, the number of associated local tensor product B-splines equals the total dimension corresponding to these T $l$-edges, as given by \eqref{eq T l-edges}.
\end{itemize}

\begin{lemma}\label{lem existence for ray}
Given a T-mesh \(\mathscr{T}\), consider its extended T-mesh with all T \(l\)-edges removed, denoted \(\mathrm{ext}_s(\mathscr{T}) \backslash T(\mathrm{ext}_s(\mathscr{T}))\). Let the set of rays in \(\mathrm{ext}_s(\mathscr{T}) \backslash T(\mathrm{ext}_s(\mathscr{T}))\) be \(\{p_1, p_2, \ldots, p_m\}\), and let \(\tilde{p}_i\) represent the ray \(p_i\) with intersections of \(p_i\) and \(p_j\) removed for \(i = 1, 2, \ldots, m\) and \(0 < j < i\). This setup satisfies the following condition:
\begin{itemize}
    \item It includes at least \(n(\tilde{p}_i)\) local tensor product B-splines associated with each \(\tilde{p}_i\) across \(n(\tilde{p}_i)\) distinct vertex sets \(V_1, \ldots, V_{n(\tilde{p}_i)}\), and for all \(1 \le j < i\), the \(n(\tilde{p}_i)\) local tensor product B-splines associated with \(\tilde{p}_i\) are not associated with \(\tilde{p}_j\). Moreover, it is possible to select appropriate \(V_1, \ldots, V_{n(\tilde{p}_i)}\) such that the one-dimensional B-splines defined by the corresponding nodes are linearly independent.
\end{itemize}
\end{lemma}

\begin{Proof}
For each ray \(\tilde{p}_i\), because each point on it is connected with another ray, one can choose the interior nodes combined with \(d_k + 1, k=1,2\) boundary nodes form \(n(\tilde{p}_i)\) local tensor product B-splines associated with \(\tilde{p}_i\). By Lemma~\ref{lem one dimension indp}, the linear independence of the one-dimensional B-splines is evident, thereby satisfying the stated condition.
\end{Proof}

As a result of Lemma \ref{lem existence} and \ref{lem existence for ray}, by extending the edges of $\mathscr{T}$, $\mathrm{ext}(\mathscr{T})$ forms a diagonalizable T-mesh, allowing the assignment of corresponding local tensor product B-splines to each edge.

It is noted that expanding a large number of edges is not always necessary to obtain an extended T-mesh that satisfies the conditions of Lemma~\ref{lem existence}. 

\begin{figure}[htbp]  
    \centering    %
    \subfigure[] %
    {
        \begin{minipage}[t]{0.33\textwidth}
            \centering      %
            \includegraphics[width=1\textwidth]{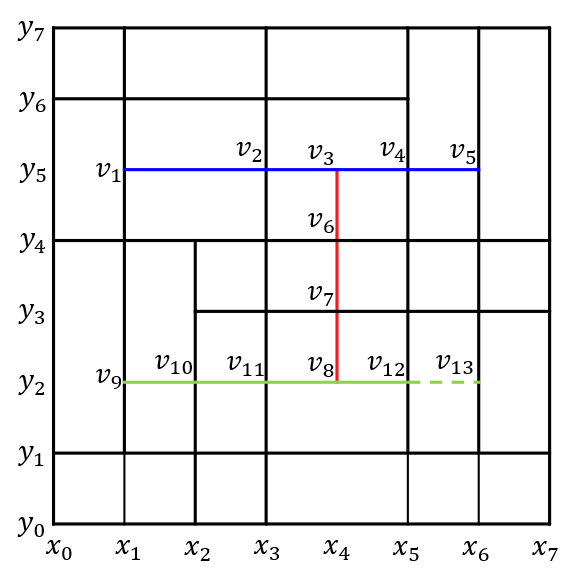}  
        \end{minipage}
    }%
    \subfigure[] %
    {
        \begin{minipage}[t]{0.33\textwidth}
            \centering      %
            \includegraphics[width=1\textwidth]{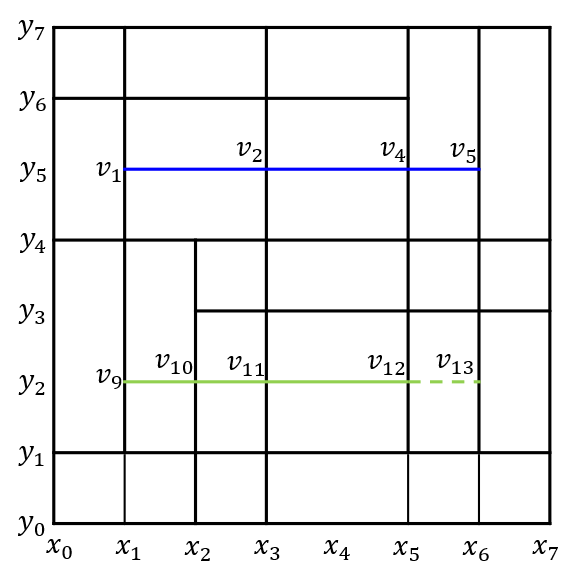}  
        \end{minipage}
    }%
    \subfigure[] %
    {
        \begin{minipage}[t]{0.33\textwidth}
            \centering      %
            \includegraphics[width=1\textwidth]{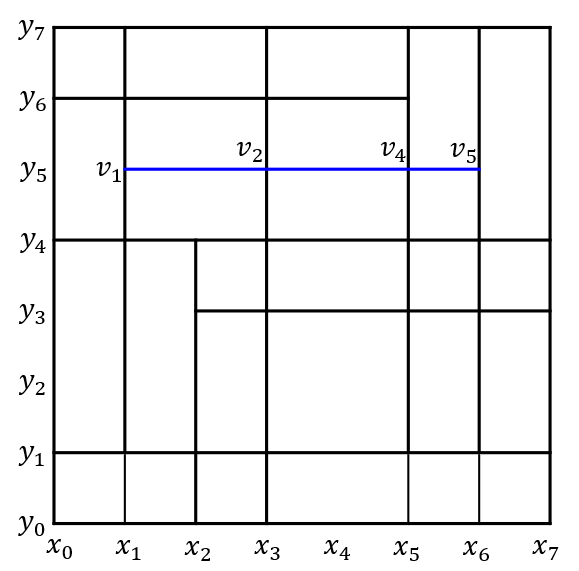}  
        \end{minipage}
    }%
    \captionsetup{font={scriptsize}}
    \caption{\textcolor{blue}{The extended T-mesh which satisfies Lemma~\ref{lem existence}. The red solid line represents $l_1$, the green solid line represents $l_2$, the blue solid line represents $l_3$, and the green dashed line denotes the temporarily extended edge. In Fig.(a), appropriate tensor-product B-splines are sought for the red T $l$-edges; in Fig.(b), appropriate tensor-product B-splines are sought for the green T $l$-edges (both solid and dashed); and in Fig.(c), appropriate tensor-product B-splines are sought for the blue T $l$-edges.}} 
    \label{wholeex3}  %
\end{figure}

\begin{example}\label{extend examp}
 We adopt the same notation as in Example \ref{tpart}. The T \(l\)-edge \(l_2\) in the T-mesh of Figure~\ref{wholeex1} lacks a local tensor product structure associated with \(l_2\) at knots \(\{v_3, v_6, v_7, v_8\}\). Consequently, we extend this T-mesh to the configuration shown in Figure~\ref{wholeex3}(a), where the blue dashed line represents the extended edge. This extended T-mesh satisfies the conditions of Lemma~\ref{lem existence}. The local tensor product B-spline associated with \(l_1\) is given by
\[
N(x_3, x_4, x_5, x_6)(x) \times N(y_2, y_3, y_4, y_5)(y).
\]
In Figure~\ref{wholeex3}(b), \(l_2\) becomes a T \(l\)-edge with endpoints \(v_9\) and \(v_{13}\), and the associated local tensor product B-splines are
\[
N(x_1, x_2, x_3, x_5)(x) \times N(y_0, y_1, y_2, y_4)(y); \quad N(x_2, x_3, x_5, x_6)(x) \times N(y_1, y_2, y_3, y_4)(y).
\]
For \(l_3\), the associated local tensor product B-spline is
\[
N(x_1, x_3, x_5, x_6)(x) \times N(y_0, y_1, y_4, y_5)(y).
\]
\end{example}

\begin{theorem}\label{thm ex linear independent}
    Given an arbitrary T-mesh $\mathscr{T}$, suppose $\mathrm{ext}_s(\mathscr{T})$ is the extended T-mesh satisfying the conditions of Lemma~\ref{lem existence}. We then consider the following three types of local tensor product B-splines:
\begin{itemize}
\item For cross-cuts in $\mathrm{ext}_s(\mathscr{T})$, we select local tensor product B-splines defined by consecutive nodes in both $x$ and $y$ directions within the tensor product part of $\mathrm{ext}_s(\mathscr{T})$, denoted as
$$\{B_i^{\mathrm{cross}}(x,y)\}_{i=1}^{\alpha},$$
where $\alpha:=(d+1+c'_h)(d+1+c'_v)$, $c'_h$ and $c'_v$ represent the numbers of horizontal and vertical cross-cuts, respectively.
\item For T $l$-edges in $\mathrm{ext}_s(\mathscr{T})$, Lemma~\ref{lem existence} allows the selection of exactly $n(\overline{l}_i) - d_k - 1$ local tensor product B-splines associated with each T $l$-edge $\overline{l}_i$ which satisfy the condition in lemma \ref{lem existence} in a $t'$-partition, for $i = 1, 2, \ldots, t'$, denoted as
$$\{B_j^T(x,y)\}_{j=1}^{\beta},$$
where $\beta:=n'_T - (d_1+1)t'_h - (d_2+1)t'_v$, $n'_T$ is the total number of vertices of all T $l$-edges in $\mathrm{ext}_s(\mathscr{T})$, $t'$ is the total number of T $l$-edges, and $t'_h$ and $t'_v$ are the numbers of horizontal and vertical T $l$-edges, respectively.

\item For rays in $\mathrm{ext}_s(\mathscr{T})\backslash T(\mathrm{ext}_s(\mathscr{T}))$. Lemma~\ref{lem existence for ray} enables the selection of exactly $n(\tilde{p}_i)$ local tensor product B-splines associated with each ray $\tilde{p}_i$,\(n(\tilde{p}_i)\) indicating the number of interior vertices on \(\tilde{p}_i\), denoted as
$$\{B_k^{\mathrm{ray}}(x,y)\}_{k=1}^{\gamma},$$
where $\gamma:=n'_v - n'_T - c'_v c'_h$, $n'_v$ is the total number of interior vertices in $\mathrm{ext}_s(\mathscr{T})$.
\end{itemize}
Thus, the set $\{B_i(x,y), B_j(x,y), B_k(x,y)\}_{i=1, j=1, k=1}^{\alpha, \beta, \gamma}$ forms a basis for $S_{d_1,d_2}(\mathrm{ext}_s(\mathscr{T}))$, expressed as
$$S_{d_1,d_2}(\mathrm{ext}_s(\mathscr{T})) = \mathrm{span}\{B_i(x,y), B_j(x,y), B_k(x,y)\}_{i=1, j=1, k=1}^{\alpha, \beta, \gamma}.$$
   
\end{theorem}
\begin{Proof}
   According to the dimension formula \eqref{dim diag},
$$\dim S_{d_1,d_2}(\mathrm{ext}_s(\mathscr{T})) = \alpha + \beta + \gamma.$$
Thus, $\dim S_{d_1,d_2}(\mathrm{ext}_s(\mathscr{T}))$ equals the number of these functions, and it remains to establish their linear independence.
Assume
\begin{equation}\label{ext basis linear independent}
    \sum_{i=1}^{\alpha} a_i B_i^{\text{cross}}(x,y) + \sum_{j=1}^{\beta} b_j B_j^T(x,y) + \sum_{k=1}^{\gamma} c_k B_k^{\text{ray}}(x,y) \equiv 0.
\end{equation}    
We need to prove that $a_i = 0$, $b_j = 0$, and $c_k = 0$ for $i = 1, 2, \ldots, \alpha$, $j = 1, 2, \ldots, \beta$, and $k = 1, 2, \ldots, \gamma$.

Without loss of generality, consider a horizontal $l$-edge $l$ in $\mathrm{ext}_s(\mathscr{T})$; the proof for vertical $l$-edges follows similarly. We apply the operators $\frac{\partial_{+}^{d_2}}{\partial y}$ and $\frac{\partial_{-}^{d_2}}{\partial y}$ to both sides of equation \eqref{ext basis linear independent}, proceeding through the following three steps:
\begin{itemize}
\item[Step 1] Assume $l$ is a T $l$-edge in the $t'$-partition of $\mathrm{ext}_s(\mathscr{T})$. Without loss of generality, let the associated local tensor product B-splines be $B_1^T(x,y), B_2^T(x,y), \ldots, B_{\beta'}^T(x,y)$. This leads to the following equations:
\begin{equation}\label{T l-edge partial +}
    \sum_{i=1}^{\alpha} a_i \frac{\partial_{+}^{d_2}}{\partial y^{d_2}} B_i^{\text{cross}}(x,y) + \sum_{j=1}^{\beta'} b_j \frac{\partial_{+}^{d_2}}{\partial y^{d_2}} B_j^T(x,y) + \sum_{j=\beta'+1}^{\beta} b_j \frac{\partial_{+}^{d_2}}{\partial y^{d_2}} B_j^T(x,y) + \sum_{k=1}^{\gamma} c_k \frac{\partial_{+}^{d_2}}{\partial y^{d_2}} B_k^{\text{ray}}(x,y) \equiv 0,
\end{equation} 

\begin{equation}\label{T l-edge partial -}
    \sum_{i=1}^{\alpha} a_i \frac{\partial_{-}^{d_2}}{\partial y^{d_2}} B_i^{\text{cross}}(x,y) + \sum_{j=1}^{\beta'} b_j \frac{\partial_{-}^{d_2}}{\partial y^{d_2}} B_j^T(x,y) + \sum_{j=\beta'+1}^{\beta} b_j \frac{\partial_{-}^{d_2}}{\partial y^{d_2}} B_j^T(x,y) + \sum_{k=1}^{\gamma} c_k \frac{\partial_{-}^{d_2}}{\partial y^{d_2}} B_k^{\text{ray}}(x,y) \equiv 0.
\end{equation}
Since every function in the equation possesses $C^{d_2-1}$ smoothness along the $y$-direction, all local tensor product B-splines except those associated with $l$ are $C^{\infty}$ across $l$. Restrict the equation on \(l\) and subtracting \eqref{T l-edge partial -} from \eqref{T l-edge partial +} yields
$$\sum_{j=1}^{\beta'} b_j \left(\frac{\partial_{+}^{d_2}}{\partial y^{d_2}} B_j^T(x,y) - \frac{\partial_{-}^{d_2}}{\partial y^{d_2}} B_j^T(x,y)\right) \equiv 0.$$
Let $B_j^T(x,y) = N_j^T(x) N_j^T(y)$, where $N_j^T(x)$ and $N_j^T(y)$ are B-spline functions in the respective directions. Then, $\frac{\partial_{+}^{d_2}}{\partial y} B_j^T(x,y) = B_j^+ N_j^T(x)$ and $\frac{\partial_{-}^{d_2}}{\partial y} B_j^T(x,y) = B_j^- N_j^T(x)$, with $B_j^+$ and $B_j^-$ being distinct constants. This reduces the equation to
$$\sum_{j=1}^{\beta'} b_j (B_j^+ - B_j^-) N_j^T(x) \equiv 0.$$
Based on the conditions of this Theorem, $\{N_j^T(x)\}_{j=1}^{\beta'}$ is a series of B-splines ensuring linear independence. It follows that $b_j (B_j^+ - B_j^-) = 0$ for $j = 1, 2, \ldots, \beta'$. Since $B_j^+ \neq B_j^-$, we conclude $b_j = 0$ for $j = 1, 2, \ldots, \beta'$.

Since local tensor product B-splines associated with distinct T $l$-edges are unique, this process applies to each T $l$-edge in the $t'$-partition of $\mathrm{ext}_s(\mathscr{T})$ following the order, yielding 
$$b_j = 0, j = 1, 2, \ldots, \beta$$

    \item[Step 2] Assume \(p\) is a ray in \(\mathrm{ext}_s(\mathscr{T}) \setminus T(\mathrm{ext}_s(\mathscr{T}))\). Without loss of generality, let the associated local tensor product B-splines be \(B_1^{\text{ray}}(x,y), B_2^{\text{ray}}(x,y), \ldots, B_{\gamma'}^{\text{ray}}(x,y)\). This leads to the following equations:
\begin{equation}\label{ray partial +}
    \sum_{i=1}^{\alpha} a_i \frac{\partial_{+}^{d_2}}{\partial y} B_i^{\text{cross}}(x,y) + \sum_{k=1}^{\gamma'} c_k \frac{\partial_{+}^{d_2}}{\partial y} B_k^{\text{ray}}(x,y) + \sum_{k=\gamma'+1}^{\gamma} c_k \frac{\partial_{+}^{d_2}}{\partial y} B_k^{\text{ray}}(x,y) \equiv 0,
\end{equation}
\begin{equation}\label{ray partial -}
    \sum_{i=1}^{\alpha} a_i \frac{\partial_{-}^{d_2}}{\partial y^{d_2}} B_i^{\text{cross}}(x,y) + \sum_{k=1}^{\gamma'} c_k \frac{\partial_{-}^{d_2}}{\partial y^{d_2}} B_k^{\text{ray}}(x,y) + \sum_{k=\gamma'+1}^{\gamma} c_k \frac{\partial_{-}^{d_2}}{\partial y^{d_2}} B_k^{\text{ray}}(x,y) \equiv 0.
\end{equation}
Since every function in this equation exhibits \(C^{d_2-1}\) smoothness along the \(y\)-direction, all local tensor product B-splines except those associated with \(p\) are \(C^{\infty}\) across \(p\). Restrict the equation on \(p\) and subtracting \eqref{ray partial -} from \eqref{ray partial +} yields
\begin{equation}
    \sum_{k=1}^{\gamma'} c_k \left(\frac{\partial_{+}^{d_2}}{\partial y^{d_2}} B_k^{\text{ray}}(x,y) - \frac{\partial_{-}^{d_2}}{\partial y^{d_2}} B_k^{\text{ray}}(x,y)\right) \equiv 0.
\end{equation}
Let \(B_k^{\text{ray}}(x,y) = N_k^{\text{ray}}(x) N_k^{\text{ray}}(y)\), where \(N_k^{\text{ray}}(x)\) and \(N_k^{\text{ray}}(y)\) are B-spline functions in the respective directions. Then, \(\frac{\partial_{+}^{d_2}}{\partial y^{d_2}} B_k^{\text{ray}}(x,y) = C_k^+ N_k^{\text{ray}}(x)\) and \(\frac{\partial_{-}^{d_2}}{\partial y^{d_2}} B_k^{\text{ray}}(x,y) = C_k^- N_k^{\text{ray}}(x)\), with \(C_k^+\) and \(C_k^-\) being distinct constants. This simplifies to
\[
\sum_{k=1}^{\gamma'} c_k (C_k^+ - C_k^-) N_k^{\text{ray}}(x) \equiv 0.
\]
Under the theorem's conditions, \(\{N_k^{\text{ray}}(x)\}_{k=1}^{\gamma'}\) consists of B-splines ensuring linear independence. Thus, \(c_k (C_k^+ - C_k^-) = 0\) for \(k = 1, 2, \ldots, \gamma'\). Since \(C_k^+ \neq C_k^-\), it follows that \(c_k = 0\) for \(k = 1, 2, \ldots, \gamma'\).

Given that local tensor product B-splines associated with distinct rays are unique, this process applies to each ray in \(\mathrm{ext}_s(\mathscr{T}) \setminus T(\mathrm{ext}_s(\mathscr{T}))\), yielding \(c_k = 0\) for \(k = 1, 2, \ldots, \gamma\).

\item[Step 3] Equation \eqref{ext basis linear independent} then reduces to
\[
\sum_{i=1}^{\alpha} a_i B_i^{\text{cross}}(x,y)\equiv 0.
\]
Due to the linear independence of \(\{B_i^{\text{cross}}(x,y)\}_{i=1}^{\alpha}\), it follows that \(a_i = 0\) for \(i = 1, 2, \ldots, \alpha\).
\end{itemize}

In conclusion, the set \(\{B_i(x,y), B_j(x,y), B_k(x,y)\}_{i=1, j=1, k=1}^{\alpha, \beta, \gamma}\) forms a basis for \(S_{d_1,d_2}(\mathrm{ext}_s(\mathscr{T}))\).
\end{Proof}

\subsection{ Extended edge elimination conditions}\label{subsection 3.2}
In Subsection~\ref{subsection 3.1}, we demonstrate that for any arbitrary T-mesh $\mathscr{T}$, there exists a diagonalizable T-mesh within its extended T-meshes, denoted as $\mathrm{ext}_s(\mathscr{T})$, such that a set of basis functions for the highest-order smooth spline space over $\mathrm{ext}_s(\mathscr{T})$ can be directly constructed from local tensor product B-splines associated with each $l$-edge. To derive a set of basis functions for the spline space over the original T-mesh $\mathscr{T}$, we must apply operations to eliminate extended edges. This subsection focuses on the method to obtain a set of spline space basis functions over $\mathscr{T}$ through these edge elimination operations.

The following lemma obviously holds.

\begin{lemma}\label{lem ascending chain}
    For the ascending chain:
    $$\mathrm{ext}_0(\mathscr{T})=\mathscr{T}\subseteq \mathrm{ext}_1(\mathscr{T})\subseteq \mathrm{ext}_2(\mathscr{T})\subseteq\cdots\subseteq\mathrm{ext}_{r}(\mathscr{T}).$$
    we have
    $$S_{d_1,d_2}\left(\mathrm{ext}_0(\mathscr{T})\right)\subseteq S_{d_1,d_2}\left(\mathrm{ext}_1(\mathscr{T})\right)\subseteq\cdots\subseteq S_{d_1,d_2}\left(\mathrm{ext}_{r}(\mathscr{T})\right).$$
    where ``$\subseteq$'' denotes a subspace relation.
\end{lemma}

In the following, we define the eliminate extended edge condition for an edge, designed to remove the smoothness constraint at that edge.

\begin{definition}
   Given a T-mesh \(\mathscr{T}\), consider the extended T-mesh \(\mathrm{ext}_s(\mathscr{T})\), where \(n = \dim S_{d_1,d_2}(\mathrm{ext}_s(\mathscr{T}))\) and \(\{\mathcal{B}_1, \mathcal{B}_2, \ldots, \mathcal{B}_n\}\) forms a basis for \(S_{d_1,d_2}(\mathrm{ext}_s(\mathscr{T}))\). For an edge \(e \in \mathrm{ext}_s(\mathscr{T})\) with \(e \notin \mathscr{T}\), we examine the following cases:
\begin{itemize}
    \item If \(e\) is a horizontal edge with equation \(y = y_k\),
    \begin{equation}\label{horizontal eee condition}
        \frac{\partial_{+}^{d_2}}{\partial y^{d_2}} \left( \sum_{i=1}^n c_i \mathcal{B}_i \right) \Big|_{e} \equiv \frac{\partial_{-}^{d_2}}{\partial y^{d_2}} \left( \sum_{i=1}^n c_i \mathcal{B}_i \right) \Big|_{e}.
    \end{equation}
    \item If \(e\) is a vertical edge with equation \(x = x_k\),
    \begin{equation}\label{vertical eee condition}
        \frac{\partial_{+}^{d_1}}{\partial x^{d_1}} \left( \sum_{i=1}^n c_i \mathcal{B}_i \right) \Big|_{e} \equiv \frac{\partial_{-}^{d_1}}{\partial x^{d_1}} \left( \sum_{i=1}^n c_i \mathcal{B}_i \right) \Big|_{e}.
    \end{equation}
\end{itemize}
We refer to the conditions~\eqref{horizontal eee condition} and~\eqref{vertical eee condition} as the \textbf{Extended Edge Elimination Condition of edge \(e\)}.
\end{definition}

For the extended edge $ l $, the eliminate extended edge condition ensures that the linear combination of $\{\mathcal{B}_1, \mathcal{B}_2, \ldots, \mathcal{B}_n\}$ achieves $C^{d_k}$ smoothness, where $k = 1$ for vertical edges and $k = 2$ for horizontal edges, given that $\mathcal{B}_i \in \mathbb{P}_{d_1,d_2}$ for $i = 1, 2, \ldots, n$, thereby rendering the linear combination $C^\infty$. In the following, we define the global eliminate extended edge condition for all extended edges.

\begin{definition}\label{EEE}
Given a T-mesh \(\mathscr{T}\), consider the extended T-mesh \(\mathrm{ext}_s(\mathscr{T})\), where \(n = \dim S_{d_1,d_2}(\mathrm{ext}_s(\mathscr{T}))\) and \(\{\mathcal{B}_1, \mathcal{B}_2, \ldots, \mathcal{B}_n\}\) forms a basis for \(S_{d_1,d_2}(\mathrm{ext}_s(\mathscr{T}))\). Let \(\{e_1, e_2, \ldots, e_s\}\) represent the set of all extended edges. For each \(e_i \in \{e_1, e_2, \ldots, e_s\}\), the eliminate extended edge conditions for \(e_i\), \(i = 1, 2, \ldots, s\), collectively form a linear system expressed as
\[
Mc \equiv 0,
\]
where \( M(i,j) = \left( \frac{\partial^{d_2}_{+}}{\partial y^{d_2}} - \frac{\partial^{d_2}_{-}}{\partial y^{d_2}} \right) \mathcal{B}_j \Big|_{e_i} \in \mathbb{P}_{d_1}[x] \), with \(\mathbb{P}_{d_1}[x]\) denoting the set of all polynomials in \(x\) of degree at most \(d_1\), if \(e_i\) is a horizontal edge, and \( M(i,j) = \left( \frac{\partial^{d_1}_{+}}{\partial x^{d_1}} - \frac{\partial^{d_1}_{-}}{\partial x^{d_1}} \right) \mathcal{B}_j \Big|_{e_i} \in \mathbb{P}_{d_2}[y] \), with \(\mathbb{P}_{d_2}[y]\) denoting the set of all polynomials in \(y\) of degree at most \(d_2\), if \(e_i\) is a vertical edge. We refer to this as the \textbf{ Extended Edge Elimination Conditions} (EEE condition for short).
\end{definition}


\begin{theorem}\label{thm eee condition basis}
     Given a T-mesh $\mathscr{T}$, consider the extended T-mesh $\mathrm{ext}_s(\mathscr{T})$, where $n = \dim S_{d_1,d_2}(\mathrm{ext}_s(\mathscr{T}))$ and $\{\mathcal{B}_1, \mathcal{B}_2, \ldots, \mathcal{B}_n\}$ forms a basis for $S_{d_1,d_2}(\mathrm{ext}_s(\mathscr{T}))$. If the homogeneous linear system corresponding to the EEE condition $Mc = 0$ admits a basis of solutions $\{\tilde{c}_j\}_{j=1}^h$, where $\tilde{c}_j = (\tilde{c}_{1,j}, \tilde{c}_{2,j}, \ldots, \tilde{c}_{n,j})$, then the functions
$$\left\{ \sum_{i=1}^n \tilde{c}_{i,j} \mathcal{B}_i \right\}_{j=1}^h.$$
form a basis for $S_{d_1,d_2}(\mathscr{T})$, and consequently, $h = \dim S_{d_1,d_2}(\mathscr{T})$.
\end{theorem}
\begin{Proof}
  Let \(\{e_1, e_2, \ldots, e_s\}\) denote the set of all extended edges. The proof proceeds in two steps:

\begin{enumerate}
    \item For any \(f \in S_{d_1,d_2}(\mathscr{T})\), Lemma~\ref{lem ascending chain} implies \(f \in S_{d_1,d_2}(\mathrm{ext}_s(\mathscr{T}))\), allowing the representation
      \[
      f = \sum_{i=1}^n f_i \mathcal{B}_i.
      \]
      Since \(f \in S_{d_1,d_2}(\mathscr{T})\), for each virtual edge \(e_i\), \(i = 1, 2, \ldots, s\), in \(\mathscr{T}\), the following conditions hold:
      \begin{itemize}
          \item If \(e_i\) is a horizontal edge with equation \(y = y_k\),
            \[
            \frac{\partial_{+}^{d_2}}{\partial y^{d_2}} \left( \sum_{i=1}^n f_i \mathcal{B}_i \right) \Big|_{y=y_k} \equiv \frac{\partial_{-}^{d_2}}{\partial y^{d_2}} \left( \sum_{i=1}^n f_i \mathcal{B}_i \right) \Big|_{y=y_k}.
            \]
          \item If \(e_i\) is a vertical edge with equation \(x = x_k\),
            \[
            \frac{\partial_{+}^{d_1}}{\partial x^{d_1}} \left( \sum_{i=1}^n f_i \mathcal{B}_i \right) \Big|_{x=x_k} \equiv \frac{\partial_{-}^{d_1}}{\partial x^{d_1}} \left( \sum_{i=1}^n f_i \mathcal{B}_i \right) \Big|_{x=x_k}.
            \]
      \end{itemize}
      Thus, \((f_1, f_2, \ldots, f_n)\) satisfies the system \(Mc = 0\), yielding
      \[
      (f_1, f_2, \ldots, f_n) = \sum_{j=1}^h k_j \tilde{c}_j,
      \]
      where \(\tilde{c}_j = (\tilde{c}_{1,j}, \tilde{c}_{2,j}, \ldots, \tilde{c}_{n,j})\) and \(k_j\), \(j = 1, 2, \ldots, h\), are constants. Consequently,
      \[
      f = \sum_{j=1}^h k_j \left( \sum_{i=1}^n \tilde{c}_{i,j} \mathcal{B}_i \right),
      \]
      implying \(f \in \mathrm{span} \left\{ \sum_{i=1}^n \tilde{c}_{i,j} \mathcal{B}_i \right\}_{j=1}^h\). Hence, \(S_{d_1,d_2}(\mathscr{T}) \subseteq \mathrm{span} \left\{ \sum_{i=1}^n \tilde{c}_{i,j} \mathcal{B}_i \right\}_{j=1}^h\).

      Conversely, under the EEE conditions, \(\sum_{i=1}^n \tilde{c}_{i,j} \mathcal{B}_i \in S_{d_1,d_2}(\mathscr{T})\) for \(j = 1, 2, \ldots, h\), so \(\mathrm{span} \left\{ \sum_{i=1}^n \tilde{c}_{i,j} \mathcal{B}_i \right\}_{j=1}^h \subseteq S_{d_1,d_2}(\mathscr{T})\). Therefore,
      \[
      S_{d_1,d_2}(\mathscr{T}) = \mathrm{span} \left\{ \sum_{i=1}^n \tilde{c}_{i,j} \mathcal{B}_i \right\}_{j=1}^h.
      \]

    \item We establish the linear independence of the functions in \(\left\{ \sum_{i=1}^n \tilde{c}_{i,j} \mathcal{B}_i \right\}_{j=1}^h\).

      Assume \(\sum_{j=1}^h d_j \left( \sum_{i=1}^n \tilde{c}_{i,j} \mathcal{B}_i \right) \equiv 0\). We prove \(d_j = 0\) for \(j = 1, 2, \ldots, h\).

      This equation simplifies to
      \[
      \sum_{i=1}^n \left( \sum_{j=1}^h d_j \tilde{c}_{i,j} \right) \mathcal{B}_i = 0,
      \]
      and since \(\{\mathcal{B}_1, \mathcal{B}_2, \ldots, \mathcal{B}_n\}\) is a basis for \(S_{d_1,d_2}(\mathrm{ext}_s(\mathscr{T}))\), it follows that
      \[
      \sum_{j=1}^h d_j \tilde{c}_{i,j} = 0, \quad i = 1, 2, \ldots, n.
      \]
      This can be expressed as
      \[
      (d_1, d_2, \ldots, d_h) \begin{pmatrix} \tilde{c}_1 \\ \tilde{c}_2 \\ \vdots \\ \tilde{c}_h \end{pmatrix} = 0.
      \]
      Given that \(\{\tilde{c}_j\}_{j=1}^h\) forms a linearly independent basis of solutions, \(d_j = 0\) for \(j = 1, 2, \ldots, h\).
\end{enumerate}

In conclusion, \(S_{d_1,d_2}(\mathscr{T}) = \mathrm{span} \left\{ \sum_{i=1}^n \tilde{c}_{i,j} \mathcal{B}_i \right\}_{j=1}^h\), and the functions in \(\left\{ \sum_{i=1}^n \tilde{c}_{i,j} \mathcal{B}_i \right\}_{j=1}^h\) are linearly independent. Thus, these functions form a basis for \(S_{d_1,d_2}(\mathscr{T})\), and consequently, \(h = \dim S_{d_1,d_2}(\mathscr{T})\).
\end{Proof}

The basis is called \textbf{PT-spline basis} of spline space $S_{d_1,d_2}(\mathscr{T})$. The construction process is outlined as follows:
\begin{itemize}
    \item[Step 1.] Convert $\mathscr{T}$ into a diagonalizable T-mesh $\mathrm{ext}_s(\mathscr{T})$ satisfying the conditions of Theorem~\ref{thm ex linear independent} through edge extension.
    \item[Step 2.] Construct a basis for $S_{d_1,d_2}(\mathrm{ext}_s(\mathscr{T}))$ using appropriate local tensor product B-splines.
    \item[Step 3.] Solve the Extended Edge Elimination (EEE) condition to obtain $\{\tilde{c}_j\}_{j=1}^{h}$.
    \item[Step 4.] By Theorem~\ref{thm eee condition basis}, the set $\left\{ \sum_{i=1}^n \tilde{c}_{i,j} \mathcal{B}_i \right\}_{j=1}^h$ forms a basis for $S_{d_1,d_2}(\mathscr{T})$.
\end{itemize}
\begin{remark}
\textcolor{blue}{The above method to construct basis functions is also valid for general spline spaces $S(d_1,d_2,\alpha,\beta,\mathscr{T})$ for any $\alpha<d_1$ and $\beta<d_2$. Specifically, in Step 2, tensor-product B-splines of degrees $(d_1,d_2)$ with continuity of order $(\alpha,\beta)$ are chosen. In Step 3, the linear combinations of tensor-product B-splines in the Extended Edge Elimination condition are required to satisfy that, along the vertical extended edges, the left- and right-hand derivatives of order 
$k_1$ are equal for 
$k_1=\alpha+1,...,d_1$, and along the horizontal extended edges, the left- and right-hand derivatives of order 
$k_2$ are equal for 
$k_2=\beta+1,...,d_2$.}
\end{remark}

An example illustrating this construction is provided below.

\begin{example}\label{}
Consider the T-mesh $\mathscr{T}$ shown in Figure~\ref{wholeex3}(a) with the spline space $S_{2,2}(\mathscr{T})$. The grid node coordinates are defined as $x_i = i$ and $y_j = j$. The PT-spline basis functions corresponding to T $l$-edges, rays, and cross-cuts are computed as follows:
\begin{itemize}
    \item \textbf{T $l$-edges}: The EEE condition for eliminating the extended edge in Example~\ref{extend examp} is
    \[
    \frac{1}{2}(6-x)^2 c_1 - (6-x)^2 c_2 \equiv 0.
    \]
    A basis for the solution space is $(c_1, c_2) = (2, 1)$. The corresponding basis for $T(\mathscr{T})$ is $\{\mathcal{B}^T_1, \mathcal{B}^T_2, \mathcal{B}^T_3\}$, where
    \begin{align*}
    \mathcal{B}^T_1 &= 2N(3,4,5,6)(x) \times N(2,3,4,5)(y) + N(2,3,5,6)(x) \times N(1,2,3,4)(y), \\
    \mathcal{B}^T_2 &= N(1,2,3,5)(x) \times N(0,1,2,4)(y), \\
    \mathcal{B}^T_3 &= N(1,3,5,6)(x) \times N(0,1,4,5)(y).
    \end{align*}
    \item \textbf{Rays}: The local tensor product B-splines are
    \[
    \{\mathcal{B}^{\mathrm{ray}}_i\}_{i=1}^3=\{N(x^1_i, \dots, x^1_{i+3})(x) \times N(6,7,7,7)(y)\}_{i=1}^3, \quad (x^1_1, \dots, x^1_6) = (0,0,0,1,3,5),
    \]
    \[
    \{\mathcal{B}^{\mathrm{ray}}_i\}_{i=4}^7=\{N(x^2_i, \dots, x^2_{i+3})(x) \times N(0,0,1,3)(y)\}_{i=1}^4, \quad (x^2_1, \dots, x^2_7) = (2,3,5,6,7,7,7),
    \]
    \[
    \{\mathcal{B}^{\mathrm{ray}}_i\}_{i=8}^9=\{N(0,0,1,2)(x) \times N(y^1_i, \dots, y^1_{i+3})(y)\}_{i=1}^2, \quad (y^1_1, \dots, y^1_5) = (0,0,0,1,4).
    \]
    \item \textbf{Tensor product part}: The local tensor product B-splines are
    \[
    \{\mathcal{B}^{\mathrm{cross}}_i\}_{i=1}^{35}=\{N(x_i, \dots, x_{i+3})(x) \times N(y_j, \dots, y_{j+3})(y)\}_{i,j=1}^{7,5},
    \]
    with $(x_1, \dots, x_{10}) = (0,0,0,1,3,5,6,7,7,7)$ and $(y_1, \dots, y_8) = (0,0,0,1,4,7,7,7)$.
\end{itemize}
The number of functions associated with T $l$-edges is 3, with rays is 9, and with tensor product components is 35. Thus, the total number of PT-spline basis functions is $35 + 9 + 3 = 47$. By Equation~\eqref{dim diag}, the dimension of $S_{2,2}(\mathscr{T})$ is
\[
\dim S_{2,2}(\mathscr{T}) = (2+1)^2 + 6 \times (2+1) + 29 - 3 \times (2+1) = 47.
\]
Hence, the number of PT-spline basis functions equals the dimension of $S_{2,2}(\mathscr{T})$.
\end{example}

\textcolor{black}{Computing the EEE condition is essential for constructing the PT-spline basis. We propose an algorithm that simplifies this process by converting global EEE condition equations into a local system. As the EEE condition applies globally to all extended edges, which may vary in selection, we minimize overlap in B-splines whose supports include cells adjacent to extended edges. This ensures independence of local equations for different extended edges, yielding a block triangular matrix for the overall EEE condition. An example demonstrates this approach.}

\begin{figure}[htbp]  
    \centering      %
    \includegraphics[width=0.4\textwidth]{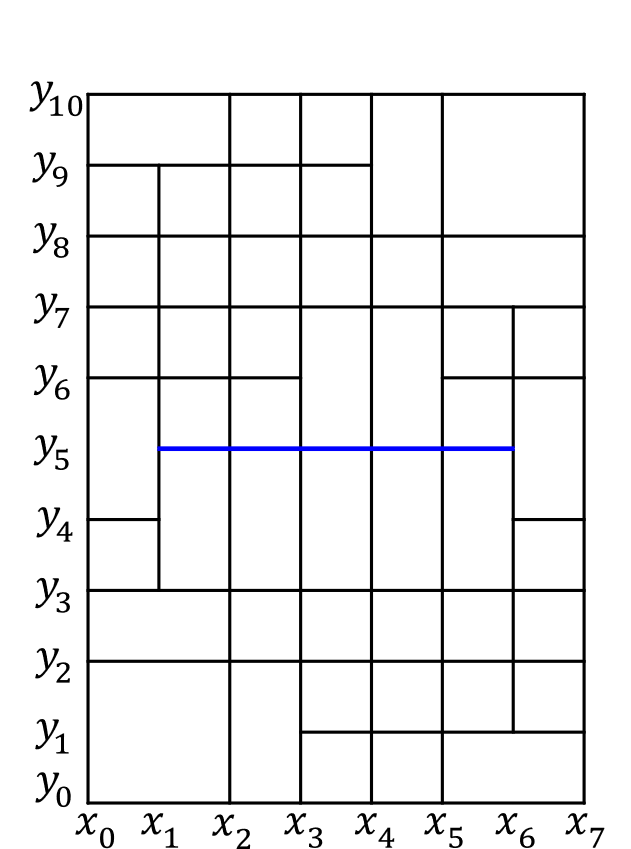} 
    \captionsetup{font={scriptsize}}
    \caption{The initial T-mesh in Example \ref{equation example}} 
    \label{equaex0}  %
\end{figure}

\begin{figure}[htbp]  
    \centering    %
    \subfigure[] %
    {
        \begin{minipage}[t]{1\textwidth}
            \centering      %
            \includegraphics[width=0.22\textwidth]{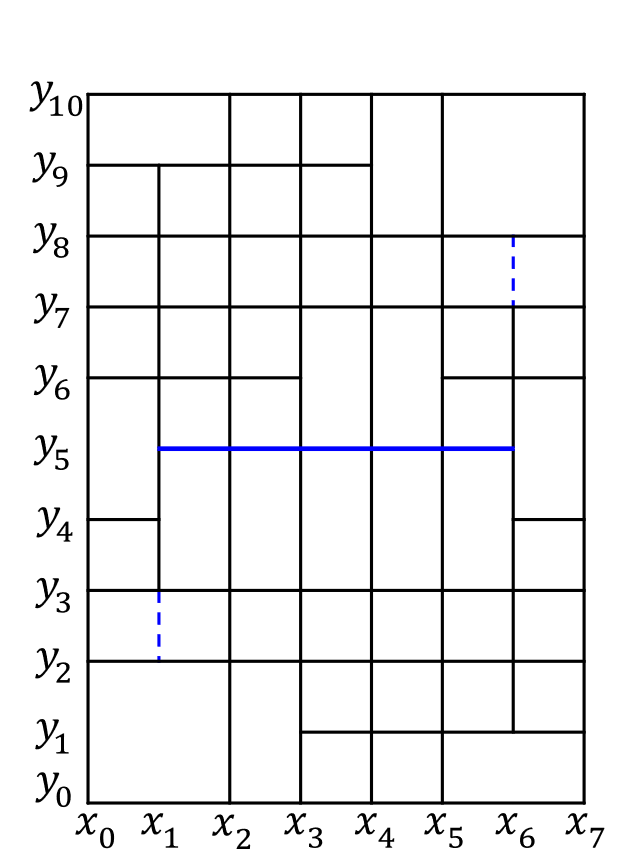}  
            \includegraphics[width=0.22\textwidth]{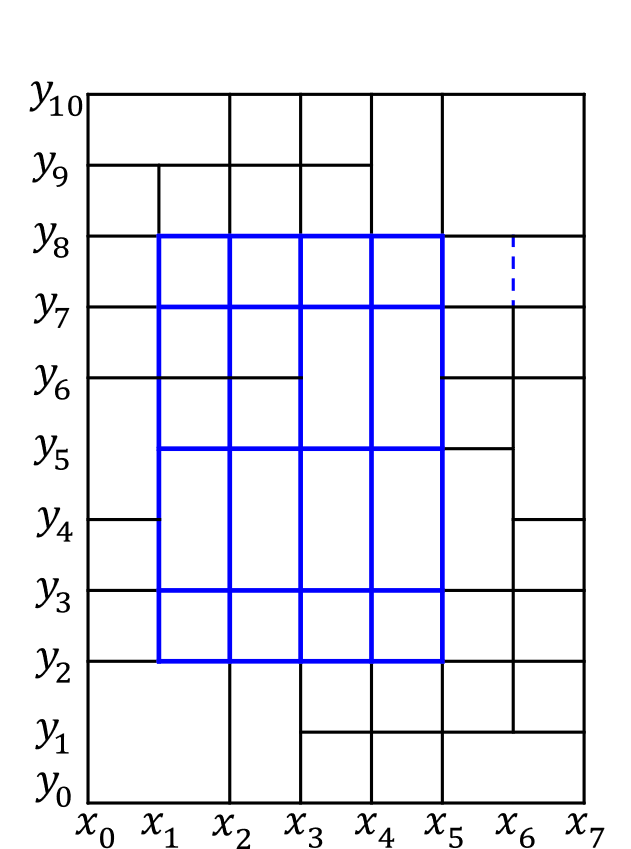} 
            \includegraphics[width=0.22\textwidth]{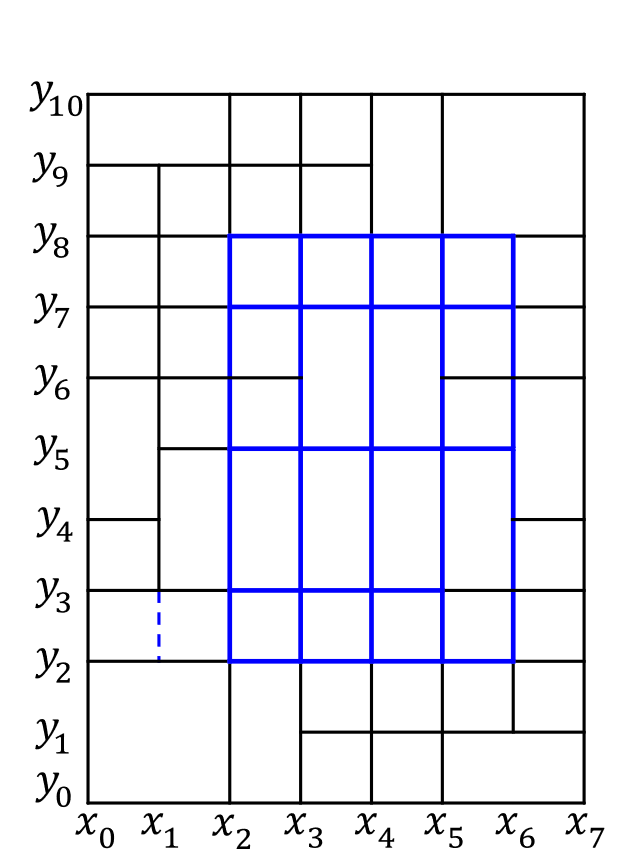} 
        \end{minipage}
    }%

    \subfigure[] %
    {
        \begin{minipage}[t]{1\textwidth}
            \centering      %
            \includegraphics[width=0.22\textwidth]{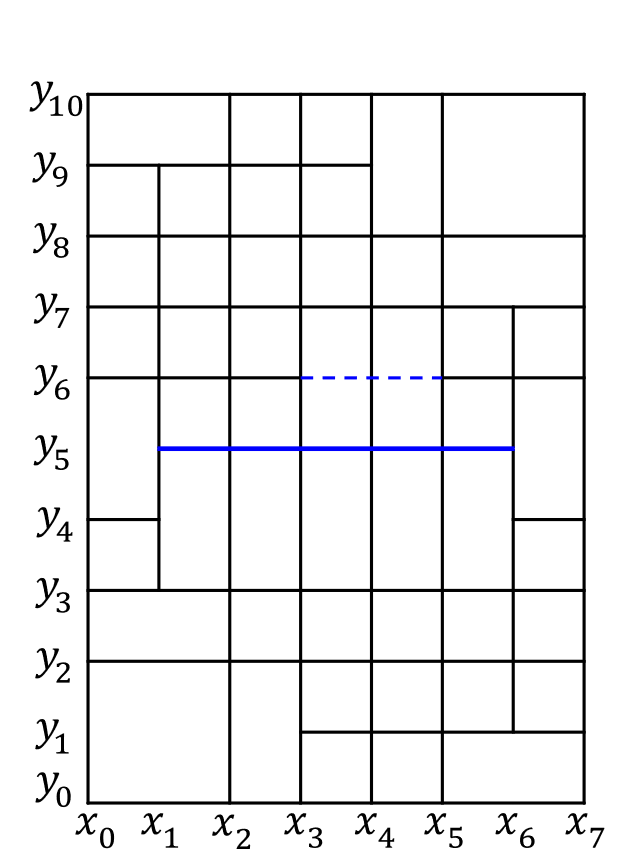}  
            \includegraphics[width=0.22\textwidth]{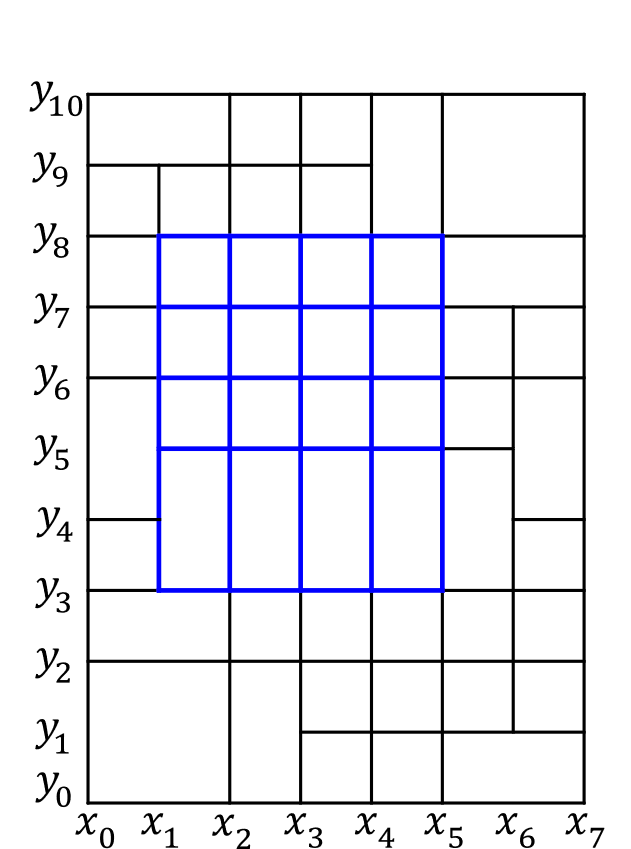} 
            \includegraphics[width=0.22\textwidth]{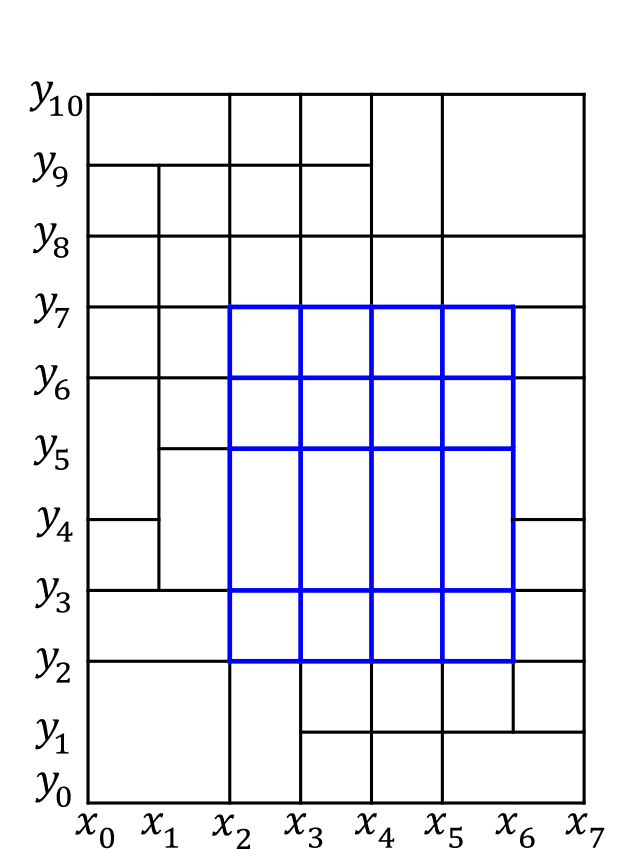} 
        \end{minipage}
    }%

    \captionsetup{font={scriptsize}}
    \caption{The extended T-meshes and associate local tensor product B-splines in Example \ref{equation example}.} 
    \label{equaex1}  %
\end{figure}

\begin{example}\label{equation example}
    Consider the T-mesh $\mathscr{T}$ depicted in Figure~\ref{equaex0} and the spline space $S_{3,3}(\mathscr{T})$. For the T $l$-edge $l$ marked in blue, the \textcolor{blue}{dimension formula~\eqref{dim diag simplify} and~\eqref{eq T l-edges}} indicates correspondence to two local tensor product B-splines. Several extensions of $\mathscr{T}$ satisfy this condition, as shown in the left subfigure of Figure~\ref{equaex1}(a) and Figure~\ref{equaex1}(b). However, the computational complexity of the EEE condition for $l$ varies significantly between these extensions. For the extension in the left subfigure of Figure~\ref{equaex1}(a), the EEE conditions for the two extended edges share no common local tensor product B-splines, allowing independent computation. Consequently, the coefficient matrix of the EEE condition’s equation system is block-diagonal locally. In contrast, the extension in Figure~\ref{equaex1}(b) may involve shared local tensor product B-splines, as illustrated in the middle and right subfigures of Figure~\ref{equaex1}(b). These shared B-splines result in EEE conditions with fewer functions while no common functions between the two edges, simplifying the equation-solving process.
\end{example}

The extension method for edges impacts the computational complexity of the EEE condition’s equation. To optimize efficiency, local tensor product B-splines in the EEE condition for extended edges should be separated as much as possible. This approach transforms the coefficient matrix into a block-diagonal form, enabling independent computation of each local extended edge and enhancing overall computational efficiency.

Finally, we establish the non-negativity of the PT-spline basis through the following proposition, ensuring all basis functions remain non-negative.

\begin{proposition}
   Given a T-mesh $\mathscr{T}$, suppose $n = \dim S_{d_1,d_2}(\mathrm{ext}_s(\mathscr{T}))$ where $\mathrm{ext}_s(\mathscr{T})$ is the extended T-mesh and $\{\mathcal{B}_1, \mathcal{B}_2, \ldots, \mathcal{B}_n\}$ forms a basis for $S_{d_1,d_2}(\mathrm{ext}_s(\mathscr{T}))$.
   and the functions
$\left\{ \sum_{i=1}^n \tilde{c}_{i,j} \mathcal{B}_i \right\}_{j=1}^h.$
form a basis for $S_{d_1,d_2}(\mathscr{T})$ where $h = \dim S_{d_1,d_2}(\mathscr{T})$.
    Then, there exist $r_{j,k}\ge 0$ for $j,k=1,2,\ldots,h$ such that 
    $$\sum\limits_{j=1}^{h}r_{j,k}\left(\sum_{i=1}^n \tilde{c}_{i,j} \mathcal{B}_i\right)\ge 0.$$
    Furthermore, 
    $$\left\{\sum\limits_{j=1}^{h}r_{j,k}\left(\sum_{i=1}^n \tilde{c}_{i,j} \mathcal{B}_i\right)\right\}_{k=1}^{h}.$$
    form a non-negative basis of the spline space $S_{d_1,d_2}(\mathscr{T})$.
\end{proposition}
\begin{Proof}
   We construct basis functions for the spline space $S_{d_1,d_2}(\mathrm{ext}_s(\mathscr{T}))$ without transforming non-cross cuts into cross-cuts. For instance, when construct associated local tensor product B-splines for a T $l$-edge needing extend some \(l\)-edges to form a cross-cut, we instead extending itself into a ray. Thus, basis functions corresponding to cross-cuts in $\mathrm{ext}_s(\mathscr{T})$, denoted $\mathcal{B}^{\mathrm{cross}}_i \in \{\mathcal{B}_1, \mathcal{B}_2, \ldots, \mathcal{B}_n\}$, $i=1,2,\ldots,\alpha$, are excluded from the Extended Edge Elimination (EEE) condition.

Assuming, without loss of generality, $\mathcal{B}^{\mathrm{cross}}_i = \mathcal{B}_i$ for $i=1,2,\ldots,\alpha$, we express
\[
\sum_{i=1}^n \tilde{c}_{i,j} \mathcal{B}_i = \sum_{i=1}^{\alpha} \mathcal{B}^{\mathrm{cross}}_i + \sum_{i=\alpha+1}^n \tilde{c}_{i,j} \mathcal{B}_i.
\]
Given $\mathcal{B}^{\mathrm{cross}}_i \geq 0$ and $\Omega \subseteq \bigcup_{i=1}^{\alpha} \mathrm{supp} \mathcal{B}^{\mathrm{cross}}_i$, where $\Omega$ is the domain of $S_{d_1,d_2}(\mathscr{T})$ for $i=1,2,\ldots,\alpha$, choosing sufficiently large $r_{j,k} \geq 0$ for $j=1,2,\ldots,h$ ensures
\[
\sum_{j=1}^{h} r_{j,k} \left( \sum_{i=1}^n \tilde{c}_{i,j} \mathcal{B}_i \right) = \sum_{j=1}^{h} \sum_{i=1}^{\alpha} r_{j,k} \mathcal{B}^{\mathrm{cross}}_i + \sum_{j=1}^{h} \sum_{i=\alpha+1}^n r_{j,k} \tilde{c}_{i,j} \mathcal{B}_i \geq 0.
\]
Additionally, the matrix $R := (r_{j,k})_{j,k=1,2,\ldots,h}$ is invertible. These conditions guarantee that the set
\[
\left\{ \sum_{j=1}^{h} r_{j,k} \left( \sum_{i=1}^n \tilde{c}_{i,j} \mathcal{B}_i \right) \right\}_{k=1}^{h}
\]
forms a non-negative basis for $S_{d_1,d_2}(\mathscr{T})$.
\end{Proof}

\section{Comparison, dimensional stability and examples}
This section presents additional examples to illustrate the benefits of PT-splines. We compare PT-spline bases with LR B-splines and provide a preliminary interpretation of dimensional instability using PT-spline bases. Further examples are provided to demonstrate the advantages of PT-spline bases over HB spline bases.

\subsection{The differences between LR B-splines and PT-splines}
In this section, we provide a comparative analysis of LR B-splines and the PT-splines. 


Locally refined B-splines (LR B-splines) extend tensor product B-splines by using local knot vectors for targeted knot insertion, enabling local mesh refinement while preserving most properties of classical B-splines~\cite{dokken2013polynomial}. However, ensuring linear independence of LR B-splines remains an unresolved challenge~\cite{dokken2013polynomial,bressan2013some,patrizi2020adaptive,patrizi2020linear}.

The primary distinction between LR B-splines and the PT-splines lies in the requirement for T-mesh expansion: PT-spline necessitate the identification of local tensor product B-splines within an extended T-mesh to construct a set of basis functions, whereas LR B-splines derive a set of local tensor product B-splines directly from the original T-mesh using a knot insertion algorithm, though these are not guaranteed to be linearly independent.

We summarize the key differences between LR B-splines and the PT-spline basis as follows:

\begin{itemize}
    \item[1.] LR B-splines may lack linear independence for a given T-mesh \(\mathscr{T}\), whereas PT-spline form a basis with the highest order of smoothness for the spline space over \(\mathscr{T}\).
    \item[2.] LR B-splines are defined solely over LR-meshes, which are diagonalizable T-meshes, representing a set of local tensor product B-splines on such meshes. In contrast, PT-spline can be defined over arbitrary T-meshes.
\end{itemize}

\begin{figure}[htbp]  
    \centering    %
    \subfigure[] %
    {
        \begin{minipage}[t]{1\textwidth}
            \centering      %
            \includegraphics[width=0.11\textwidth]{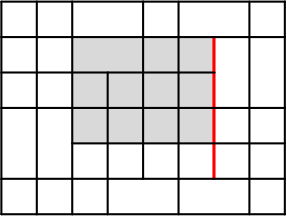}
            \includegraphics[width=0.11\textwidth]{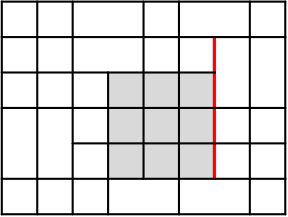}
            \includegraphics[width=0.11\textwidth]{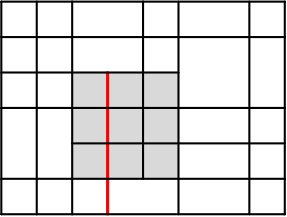}
            \includegraphics[width=0.11\textwidth]{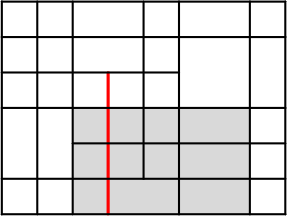}
            \includegraphics[width=0.11\textwidth]{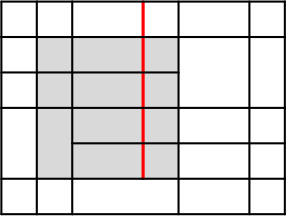}
            \includegraphics[width=0.11\textwidth]{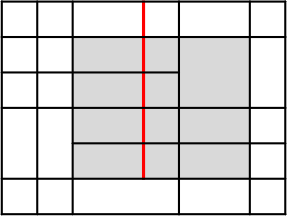}
            \includegraphics[width=0.11\textwidth]{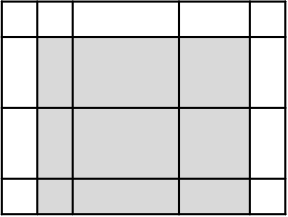}
        \end{minipage}
    }%
    
    \subfigure[] %
    {
        \begin{minipage}[t]{1\textwidth}
            \centering      %
            \includegraphics[width=0.11\textwidth]{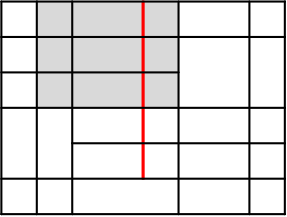} 
            \includegraphics[width=0.11\textwidth]{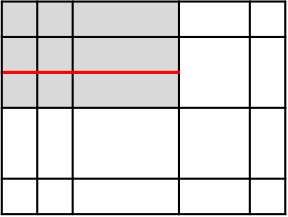} 
            \includegraphics[width=0.11\textwidth]{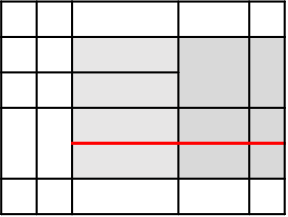} 
            \includegraphics[width=0.11\textwidth]{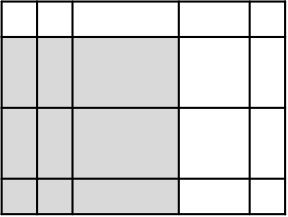}
            \includegraphics[width=0.11\textwidth]{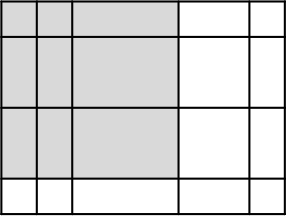} 
            \includegraphics[width=0.11\textwidth]{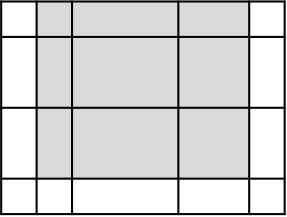}
            \includegraphics[width=0.11\textwidth]{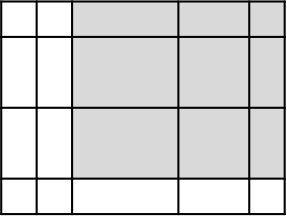}
            \includegraphics[width=0.11\textwidth]{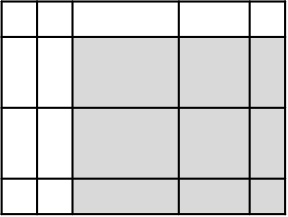} 
        \end{minipage}
    }%
    
    \captionsetup{font={scriptsize}}
    \caption{The support of PT-spline basis over T-mesh.} 
    \label{contrLRex1}  %
\end{figure}

\begin{figure}[htbp]  
    \centering    %
    \subfigure[] %
    {
        \begin{minipage}[t]{1\textwidth}
            \centering      %
            \includegraphics[width=0.11\textwidth]{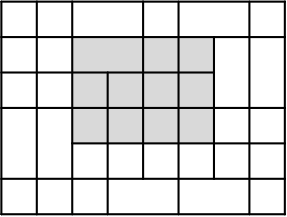}
            \includegraphics[width=0.11\textwidth]{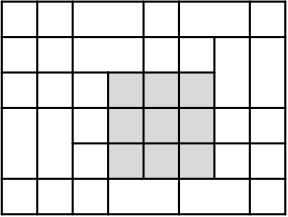}
            \includegraphics[width=0.11\textwidth]{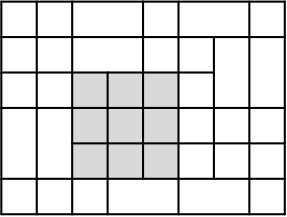}
            \includegraphics[width=0.11\textwidth]{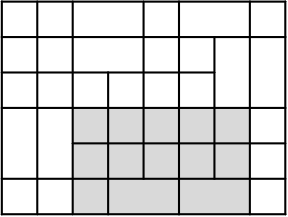}
            \includegraphics[width=0.11\textwidth]{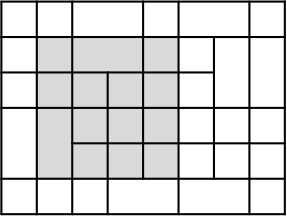}
            \includegraphics[width=0.11\textwidth]{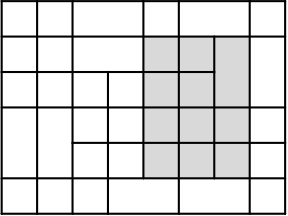}
            \includegraphics[width=0.11\textwidth]{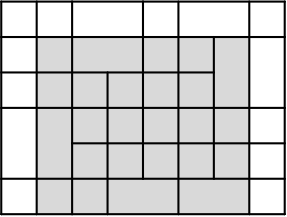}
            \includegraphics[width=0.11\textwidth]{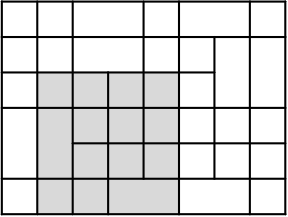} 
        \end{minipage}
    }%
    
    \subfigure[] %
    {
        \begin{minipage}[t]{1\textwidth}
            \centering      %
            \includegraphics[width=0.11\textwidth]{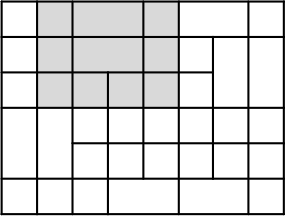} 
            \includegraphics[width=0.11\textwidth]{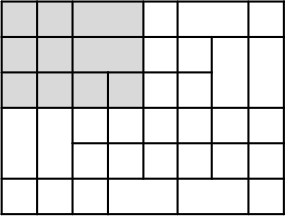} 
            \includegraphics[width=0.11\textwidth]{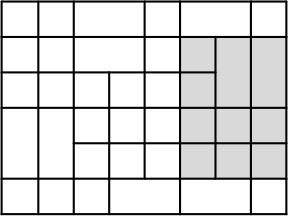}
            \includegraphics[width=0.11\textwidth]{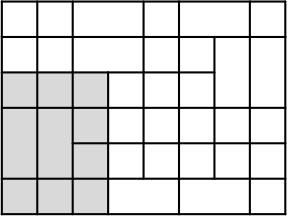} 
            \includegraphics[width=0.11\textwidth]{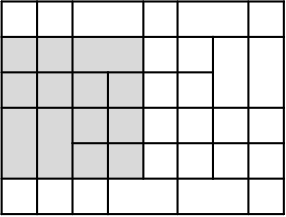}
            \includegraphics[width=0.11\textwidth]{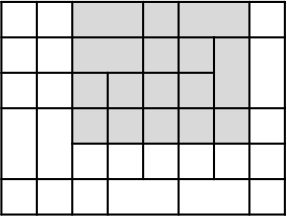}
            \includegraphics[width=0.11\textwidth]{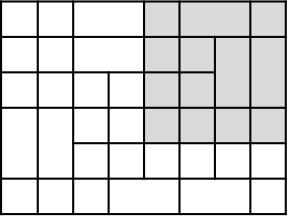} 
            \includegraphics[width=0.11\textwidth]{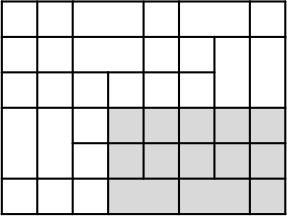} 
        \end{minipage}
    }%
    
    \captionsetup{font={scriptsize}}
    \caption{The support of LR splines over T-mesh.} 
    \label{contrLRex2}  %
\end{figure}

Below, we provide a detailed explanation of these differences, supported by examples.
\begin{example}
In this example, we introduce a T-mesh, as outlined in~\cite{dokken2013polynomial}, to compare the PT-splines with LR-splines, highlighting the linear independence of the PT-spline basis. Notably, for this T-mesh, the PT-spline basis can be constructed directly using local tensor product B-splines without extending the T-mesh, as illustrated in Figure~\ref{contrLRex1}.

In Figure~\ref{contrLRex1}, we mark the support of the PT-spline basis for the T-mesh in \textcolor{blue}{Figure~\ref{contrLRex1} (a)} with gray and indicate the lines used for spline construction with red lines, following the PT-spline basis construction procedure. As reported in~\cite{dokken2013polynomial}, LR B-spline demonstrates linear dependence, as depicted in Figure~\ref{contrLRex2} (a). In contrast, the PT-spline basis consistently ensures linear independence. \textcolor{blue}{For this particular example, it can been observed that the linear dependence of the LR B-spline collection can be resolved by replacing the B-splines in second and third picture in Figure~\ref{contrLRex2} (a) with their sum. Thus, removing linear dependence and at the same time keeping partition of unity.}
\end{example}

\subsection{Dimensional instability of spline spaces over T-meshes }

In this section, we utilize an example to elucidate the relationship between dimensional instability of spline spaces and PT-spline bases, offering a new perspective on this phenomenon.

\begin{figure}[htbp]  
    \centering    %
    \subfigure[$\mathscr{T}$] %
    {
        \begin{minipage}[t]{0.2\textwidth}
            \centering      %
            \includegraphics[width=1\textwidth]{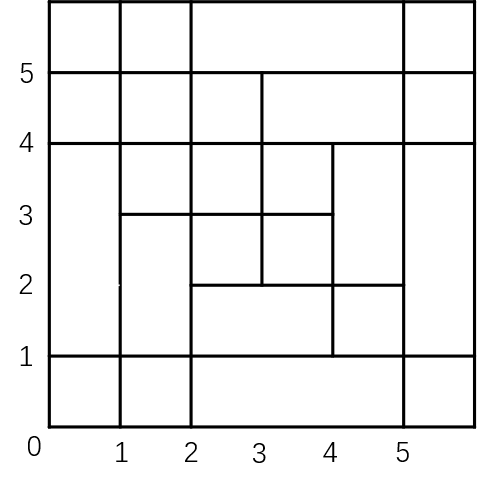}  
        \end{minipage}
    }%
    \subfigure[$\text{Support of }B_1$] %
    {
        \begin{minipage}[t]{0.2\textwidth}
            \centering      %
            \includegraphics[width=1\textwidth]{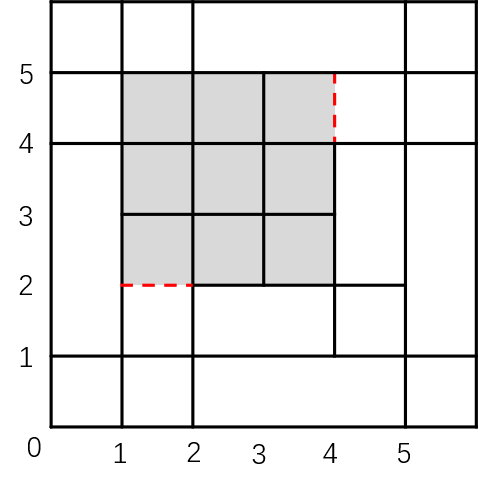}  
        \end{minipage}
    }%
    \subfigure[$\text{Support of }B_2$] %
    {
        \begin{minipage}[t]{0.2\textwidth}
            \centering      %
            \includegraphics[width=1\textwidth]{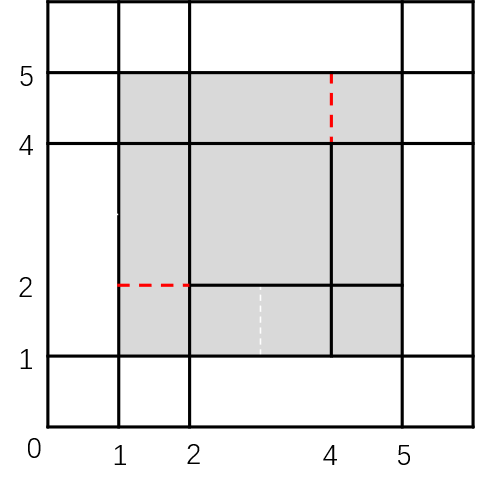}  
        \end{minipage}
    }%
    \subfigure[Tensor product mesh] %
    {
        \begin{minipage}[t]{0.2\textwidth}
            \centering      %
            \includegraphics[width=1\textwidth]{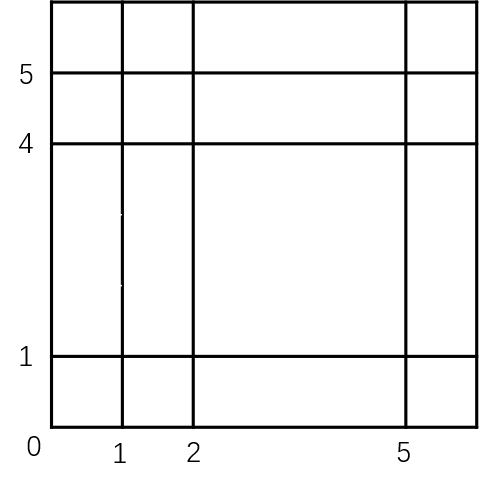}  
        \end{minipage}
    }%
    \captionsetup{font={scriptsize}}
    \caption{The process to construct PT-spline basis in Example \ref{exm ins}.} 
    \label{instaex1}  %
\end{figure}

\begin{figure}[htbp]  
    \centering      %
    \includegraphics[width=0.35\textwidth]{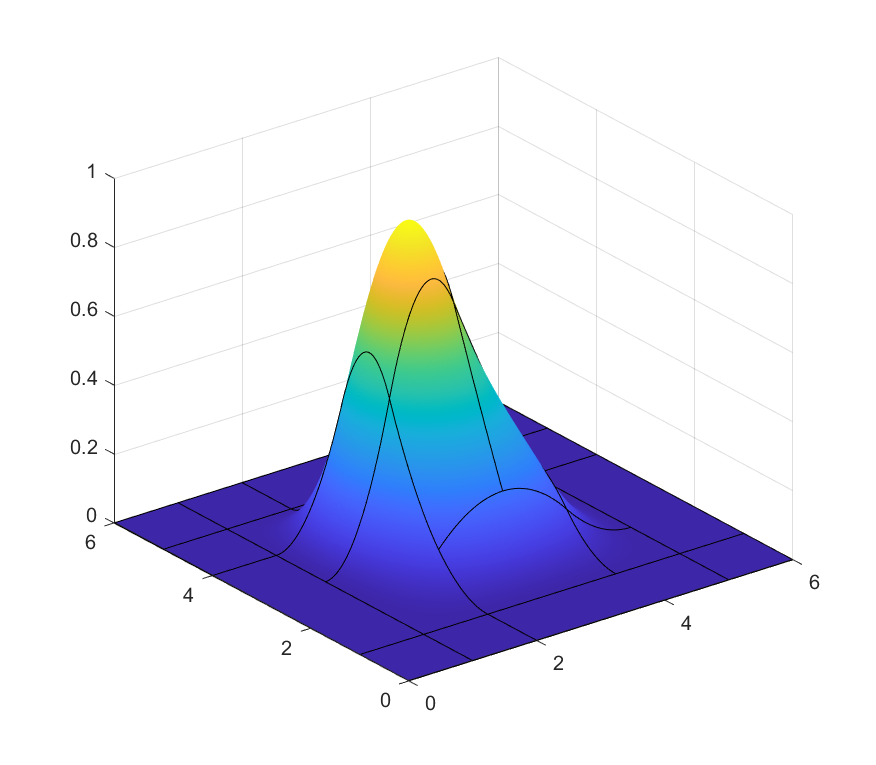} 
    \captionsetup{font={scriptsize}}
    \caption{The PT-spline basis $\frac{8}{9}B_1+B_2$.} 
    \label{Cspline1}  %
\end{figure}

\begin{figure}[htbp]  
    \centering    %
    \subfigure[$\mathscr{T}$] %
    {
        \begin{minipage}[t]{0.2\textwidth}
            \centering      %
            \includegraphics[width=1\textwidth]{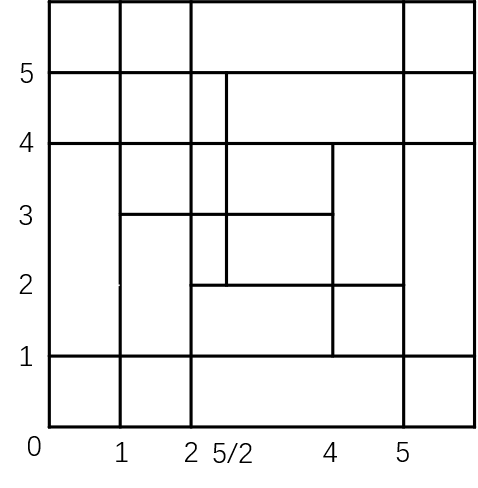}  
        \end{minipage}
    }%
    \subfigure[$\text{Support of }B_1$] %
    {
        \begin{minipage}[t]{0.2\textwidth}
            \centering      %
            \includegraphics[width=1\textwidth]{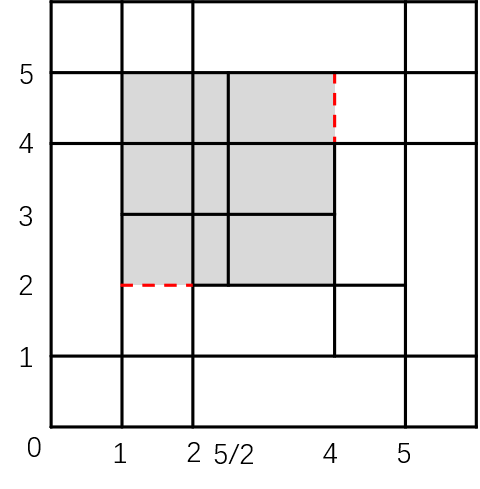}  
        \end{minipage}
    }%
    \subfigure[$\text{Support of }B_2$] %
    {
        \begin{minipage}[t]{0.2\textwidth}
            \centering      %
            \includegraphics[width=1\textwidth]{stable2c.png}  
        \end{minipage}
    }%
    \subfigure[Tensor product mesh] %
    {
        \begin{minipage}[t]{0.2\textwidth}
            \centering      %
            \includegraphics[width=1\textwidth]{stable2d.png}  
        \end{minipage}
    }%
    \captionsetup{font={scriptsize}}
    \caption{The process to construct PT-spline basis in Example~\ref{exm ins}.} 
    \label{instaex2}  %
\end{figure}

\begin{example}\label{exm ins}
Consider the T-mesh \(\mathscr{T}\) illustrated in Figure~\ref{instaex1}(a), sourced from~\cite{li2011instability}. In~\cite{li2011instability}, it is shown that the dimension of the spline space \(S_2(\mathscr{T})\) is unstable. To demonstrate this dimensional instability using the construction of the PT-spline basis, we introduce two T-meshes sharing the same topological structure as Figure~\ref{instaex1}(a), with exact coordinate data presented in Figure~\ref{instaex1}(a) and \textcolor{blue}{Figure~\ref{instaex2}(a)}. We construct the basis for these T-meshes independently.

In Figure~\ref{instaex1}, only the two functions in \textcolor{blue}{Figure~\ref{instaex1} (b) and Figure~\ref{instaex1} (c)} lack \(C^2\) continuity in the \(y\)-direction at \([1, 2] \times \{2\}\) and in the \(x\)-direction at \(\{4\} \times [4, 5]\). The EEE condition is
$$\begin{cases}
    c_1(\frac{\partial^2_+ }{\partial y^2}-\frac{\partial^2_- }{\partial y^2})B_1\Big|_{[1,2]\times \{2\}}+c_2(\frac{\partial^2_+ }{\partial y^2}-\frac{\partial^2_- }{\partial y^2})B_2\Big|_{[1,2]\times \{2\}}=0\\
    c_1(\frac{\partial^2_+ }{\partial x^2}-\frac{\partial^2_- }{\partial x^2})B_1\Big|_{\{4\}\times [4,5]}+c_2(\frac{\partial^2_+ }{\partial x^2}-\frac{\partial^2_- }{\partial x^2})B_2\Big|_{\{4\}\times [4,5]}=0
\end{cases}$$
According to Table~\ref{instable1}, that is 
$$\begin{cases}
\frac{1}{2}c_1-\frac{4}{9}c_2=0\\
\textcolor{blue}{\frac{-1}{2}c_1+\frac{4}{9}c_2=0}
\end{cases}$$
The basis for this T-mesh is $\frac{8}{9}B_1+B_2$ which is illustrated in Figure~\ref{Cspline1} along with the basis for Figure\ref{instaex1} (d) and the dimension of the T-mesh is 37.

In Figure~\ref{instaex2}, only two functions in \textcolor{blue}{Figure~\ref{instaex2} (b) and Figure~\ref{instaex2} (c)} are not $C^2$ continuity in $y$-direction at $[1,2]\times \{2\}$ and $C^2$ continuity in $x$-direction at $\{4\}\times [4,5]$. The EEE condition is 
$$\begin{cases}
    c_1(\frac{\partial^2_+ }{\partial y^2}-\frac{\partial^2_- }{\partial y^2})B_1\Big|_{[1,2]\times \{2\}}+c_2(\frac{\partial^2_+ }{\partial y^2}-\frac{\partial^2_- }{\partial y^2})B_2\Big|_{[1,2]\times \{2\}}=0\\
    c_1(\frac{\partial^2_+ }{\partial x^2}-\frac{\partial^2_- }{\partial x^2})B_1\Big|_{\{4\}\times [4,5]}+c_2(\frac{\partial^2_+ }{\partial x^2}-\frac{\partial^2_- }{\partial x^2})B_2\Big|_{\{4\}\times [4,5]}=0
\end{cases}$$
According to Table~\ref{instable2}, that is 
$$\begin{cases}
\frac{2}{3}c_1-\frac{4}{9}c_2=0\\
\frac{-1}{3}c_1+\frac{4}{9}c_2=0
\end{cases}$$
The solution space is empty, so the basis for this T-mesh is the basis in Figure~\ref{instaex2} (d) and the dimension of the T-mesh is 36.
\end{example}
\begin{table}[htbp]
\centering
\caption{Functions of $\frac{\partial^2_+ }{\partial y^2}$ and $\frac{\partial^2_- }{\partial y^2}$ for $B_1$ and $B_2$ in Example~\ref{exm ins}}
\setlength{\tabcolsep}{15pt}
\renewcommand{\arraystretch}{2}
\begin{tabular}{ccc}
\toprule
         & $B_1$          & $B_2$    \\
\bottomrule
$\frac{\partial^2_- }{\partial y^2}$         
 & 0$\times \frac{1}{2}(x-1)^2$           
 & $\frac{2}{3}\times \frac{1}{3}(x-1)^2$   \\
 $\frac{\partial^2_+ }{\partial y^2}$        
 & 1$\times \frac{1}{2}(x-1)^2$          
 & $\frac{-2}{3}\times \frac{1}{3}(x-1)^2$  \\
 $\frac{\partial^2_- }{\partial x^2}$        
 & 1$\times \frac{1}{2}(5-y)^2$         
 & $\frac{-2}{3}\times \frac{1}{3}(5-y)^2$   \\
 $\frac{\partial^2_+ }{\partial x^2}$        
 & 0$\times \frac{1}{2}(5-y)^2$         
 & $\frac{2}{3}\times \frac{1}{3}(5-y)^2$  \\
\bottomrule
\end{tabular}
\label{instable1}
\end{table}

\begin{table}[htbp]
\centering
\caption{Functions of $\frac{\partial^2_+ }{\partial y^2}$ and $\frac{\partial^2_- }{\partial y^2}$ for $B_1$ and $B_2$ in Example~\ref{exm ins}}
\setlength{\tabcolsep}{15pt}
\renewcommand{\arraystretch}{2}
\begin{tabular}{ccc}
\toprule
         & $B_1$          & $B_2$    \\
\bottomrule
$\frac{\partial^2_- }{\partial y^2}$         
 & 0$\times \frac{2}{3}(x-1)^2$           
 & $\frac{2}{3}\times \frac{1}{3}(x-1)^2$   \\
 $\frac{\partial^2_+ }{\partial y^2}$        
 & 1$\times \frac{2}{3}(x-1)^2$          
 & $\frac{-2}{3}\times \frac{1}{3}(x-1)^2$  \\
 $\frac{\partial^2_- }{\partial x^2}$        
 & $\frac{2}{3}\times \frac{1}{2}(5-y)^2$         
 & $\frac{-2}{3}\times \frac{1}{3}(5-y)^2$   \\
 $\frac{\partial^2_+ }{\partial x^2}$        
 & 0$\times \frac{1}{2}(5-y)^2$         
 & $\frac{2}{3}\times \frac{1}{3}(5-y)^2$  \\
\bottomrule
\end{tabular}
\label{instable2}
\end{table}

From Example~\ref{exm ins}, it is evident that the dimensional instability of spline spaces over T-meshes can be attributed to variations in the number of basis solutions of the EEE condition for T-meshes with identical topological structures but differing coordinate data.

\subsection{More examples of PT-spline basis construction}
In this section, we present additional examples to demonstrate the advantages of the PT-spline basis. As a first example, we show that the PT-spline basis can be constructed over any T-mesh containing T-cycles~\cite{li2011instability,guo2015problem}, marking the novel instance of basis construction on such T-meshes.

\begin{figure}[htbp]  
    \centering    %
    \subfigure[\textcolor{blue}{$\mathscr{T}$}] %
    {
        \begin{minipage}[t]{0.16\textwidth}
            \centering      %
            \includegraphics[width=1\textwidth]{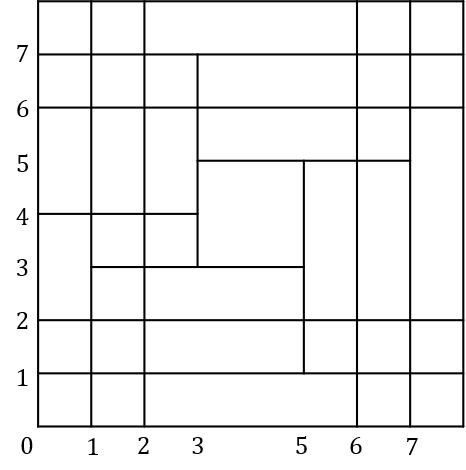}  
        \end{minipage}
    }%
    \subfigure[] %
    {
        \begin{minipage}[t]{0.16\textwidth}
            \centering      %
            \includegraphics[width=1\textwidth]{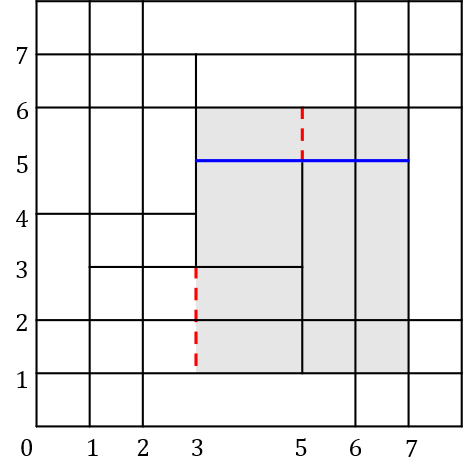}  
        \end{minipage}
    }%
    \subfigure[] %
    {
        \begin{minipage}[t]{0.16\textwidth}
            \centering      %
            \includegraphics[width=1\textwidth]{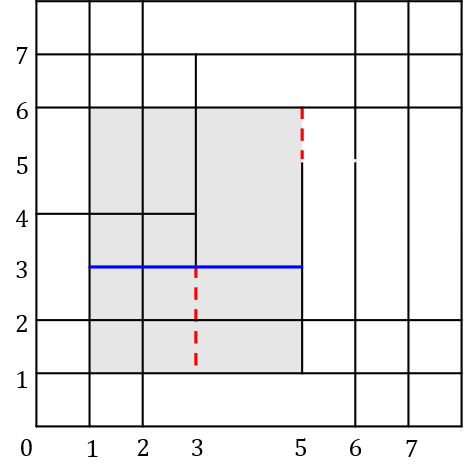}  
        \end{minipage}
    }%
    \subfigure[] %
    {
        \begin{minipage}[t]{0.16\textwidth}
            \centering      %
            \includegraphics[width=1\textwidth]{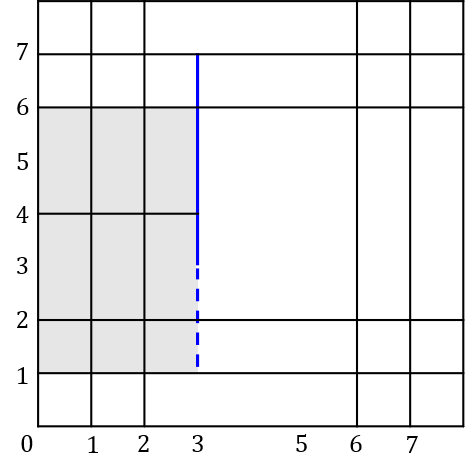}  
        \end{minipage}
    }%
    \subfigure[] %
    {
        \begin{minipage}[t]{0.16\textwidth}
            \centering      %
            \includegraphics[width=1\textwidth]{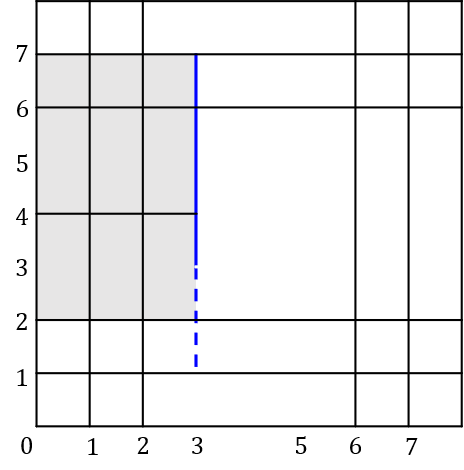}  
        \end{minipage}
    }%
    \subfigure[\textcolor{blue}{$\mathscr{T}_1$}] %
    {
        \begin{minipage}[t]{0.16\textwidth}
            \centering      %
            \includegraphics[width=1\textwidth]{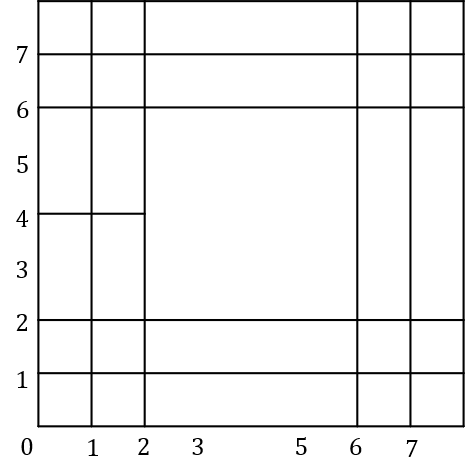}  
        \end{minipage}
    }%
    \captionsetup{font={scriptsize}}
    \caption{The example to construct PT-spline basis with T-cycle.} 
    \label{T-cycle}
\end{figure}

\begin{example}
PT-spline basis can accommodate T-meshes with T-cycles, whereas LR B-splines cannot be defined over T-meshes containing T-cycles~\cite{guo2015problem}, and PT-spline enable basis construction for such meshes. The initial T-mesh is depicted in Figure~\ref{T-cycle}(a). For the four tensor product (TP) B-splines with knot vectors
\[
B_1: [3, 5, 6, 7] \times [1, 2, 5, 6], \quad B_2: [1, 2, 3, 5] \times [1, 2, 3, 6],
\]
\[
B_3: [0, 1, 2, 3] \times [1, 2, 4, 6], \quad B_4: [0, 1, 2, 3] \times [2, 4, 6, 7],
\]
the following system of equations, derived from constraints on the red dashed lines as listed in Table~\ref{T-cycle table}, is obtained:
\[
\begin{cases}
\frac{1}{2} c_1 + \frac{1}{36} c_2 = 0, \\
\frac{1}{12} c_1 + c_2 - \frac{1}{3} c_3 = 0, \\
\frac{-1}{18} c_1 - \frac{3}{2} c_2 + \frac{7}{24} c_3 - \frac{1}{8} c_4 = 0, \\
\frac{7}{18} c_1 + 8 c_2 - \frac{11}{6} c_3 + \frac{1}{2} c_4 = 0, \\
\frac{-17}{36} c_1 - 9 c_2 + \frac{13}{6} c_3 - \frac{1}{2} c_4 = 0.
\end{cases}
\]
The solution space has a dimension of one, yielding a single function \( f = \frac{130}{651} B_1 - \frac{7}{651} B_2 + \frac{387}{649} B_3 + B_4 \) for the T-mesh in \textcolor{blue}{Figure~\ref{T-cycle} (a)}. Since \( f \) includes negative values and the function \( B_0 \) with knot vector \([1, 2, 6, 7] \times [1, 2, 6, 7]\) serves as a PT-spline basis for \( \mathscr{T}_1 \), the final PT-spline basis for \( \mathscr{T} \) comprises \( f + B_0 \) and the PT-spline basis for T-mesh \( \mathscr{T}_1 \).

\end{example}
\begin{table}[htbp]
\centering
\caption{Functions of $\frac{\partial^2_+ }{\partial x^2}$ and $\frac{\partial^2_- }{\partial x^2}$ for the four splines on three line segments}
\setlength{\tabcolsep}{15pt}
\renewcommand{\arraystretch}{2}
\begin{tabular}{ccccc}
\toprule
    &   & $5\times [5,6]$    & $3\times [1,2]$  & $3\times [2,3]$  \\
\bottomrule
\multirow{2}{*}{$B_1$} & $\frac{\partial^2_- }{\partial x^2}$         
 & $\frac{1}{3}\times \frac{(6-y)^2}{4}$ 
 & $0 \times \frac{(y-1)^2}{4}$           
 & $0 \times (\frac{-3}{4}y^2+4y-\frac{9}{2})$   \\
 & $\frac{\partial^2_+ }{\partial x^2}$        
 & $\frac{-5}{3}\times \frac{(6-y)^2}{4}$ 
 & $\frac{1}{3} \times \frac{(y-1)^2}{4}$          
 & $\frac{1}{3}\times (\frac{-3}{4}y^2+4y-\frac{9}{2})$  \\
\multirow{2}{*}{$B_2$} & $\frac{\partial^2_- }{\partial x^2}$         
 & $\frac{1}{3}\times \frac{(6-y)^2}{12}$ 
 & $\frac{-5}{3} \times \frac{(y-1)^2}{2}$           
 & $\frac{-5}{3}\times (\frac{-1}{6}y^2+\frac{7}{6}y-\frac{17}{12})$ \\
 & $\frac{\partial^2_+ }{\partial x^2}$        
 & $0\times \frac{(6-y)^2}{12}$ 
 & $\frac{1}{3} \times \frac{(y-1)^2}{2}$    
 & $\frac{1}{3} \times (\frac{-1}{6}y^2+\frac{7}{6}y-\frac{17}{12})$ \\
 \multirow{2}{*}{$B_3$} & $\frac{\partial^2_- }{\partial x^2}$         
 &  
 & $1 \times \frac{(y-1)^2}{3}$           
 & $1 \times (\frac{-7}{24}y^2+\frac{11}{6}y-\frac{13}{6})$   \\
 & $\frac{\partial^2_+ }{\partial x^2}$        
 &  
 & $0 \times \frac{(x-1)^2}{3}$          
 & $0 \times (\frac{-7}{24}y^2+\frac{11}{6}y-\frac{13}{6})$  \\
 \multirow{2}{*}{$B_4$} & $\frac{\partial^2_- }{\partial x^2}$         
 &  &            
 & $1\times (\frac{1}{8}y^2-\frac{1}{2}y+\frac{1}{2})$   \\
 & $\frac{\partial^2_+ }{\partial x^2}$        
 &  &           
 & $0\times (\frac{1}{8}y^2-\frac{1}{2}y+\frac{1}{2})$  \\
\bottomrule
\end{tabular}
\label{T-cycle table}
\end{table}

\textcolor{blue}{In example \ref{exm hier}}, we demonstrate that for certain hierarchical T-meshes, (T)HB-splines~\cite{vuong2011hierarchical,giannelli2012thb} do not form a basis for the spline space, whereas the PT-spline basis does.

\begin{figure}[htbp]  
    \centering    %
    \subfigure[] %
    {
        \begin{minipage}[t]{0.48\textwidth}
            \centering      %
            \includegraphics[width=1\textwidth]{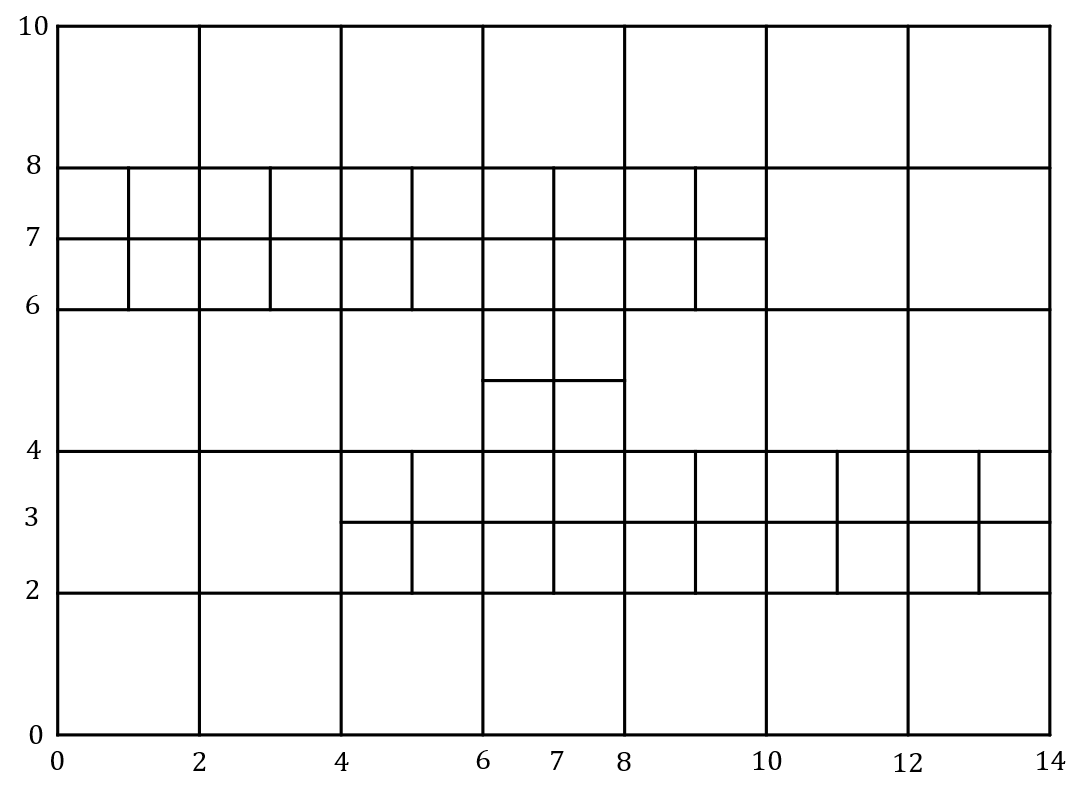}  
        \end{minipage}
    }%
    \subfigure[] %
    {
        \begin{minipage}[t]{0.48\textwidth}
            \centering      %
            \includegraphics[width=1\textwidth]{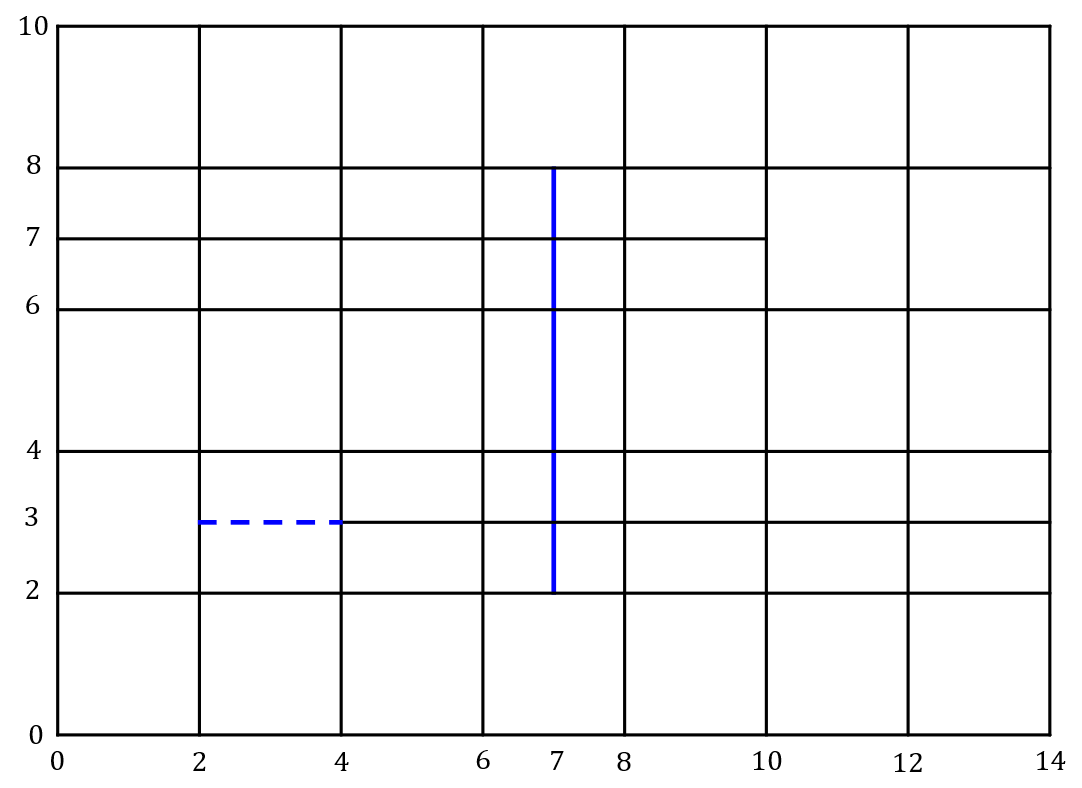}  
        \end{minipage}
    }%
    \captionsetup{font={scriptsize}}
    \caption{The hierarchical T-mesh and its extended mesh.} 
    \label{hierexfig}  %
\end{figure}

\begin{example}\label{exm hier}
Consider the hierarchical T-mesh $\mathscr{T}$ shown in Figure~\ref{hierexfig} (a) with the spline space $S_{4,4}(\mathscr{T})$. Existing hierarchical spline constructions, such as (T)HB~\cite{vuong2011hierarchical,giannelli2012thb}, generate only 99 functions at the first level, which is less than the dimension of $S_{4,4}(\mathscr{T})$, calculated as 110. We construct a basis for $S_{4,4}(\mathscr{T})$ using PT-spline. Figure~\ref{hierexfig} (b) shows the extended T-mesh of $\mathscr{T}$ after removing vanished $l$-edges. The PT-spline basis is categorized into three types:
\begin{itemize}
    \item \textbf{T $l$-edges}: The two local tensor product B-splines at the extended edge are
    \[
    N(2,4,6,7,8,10)(x) \times N(2,3,4,6,7,8)(y), \quad N(2,4,6,8,10,12)(x) \times N(0,0,0,2,3,4)(y).
    \]
    The EEE condition yields
    \[
    -\frac{9}{200} \left( \frac{x-2}{2} \right)^4 c_1 + \frac{317}{864} \left( \frac{x-2}{2} \right)^4 c_2 \equiv 0.
    \]
    One basis of solution space is $(c_1,c_2)=(1, \frac{972}{7925})$, thus, one PT-spline is
    \[
    \mathcal{B}_1^T = N(2,4,6,7,8,10)(x) \times N(2,3,4,6,7,8)(y) + \frac{972}{7925} N(2,4,6,8,10,12)(x) \times N(0,0,0,2,3,4)(y).
    \]
    \item \textbf{Rays}: The tensor product B-splines are
    \[
    \{\mathcal{B}^{\mathrm{ray}}_i\}_{i=1}^{5} = \{ N(x^1_i, \dots, x^1_{i+5})(x) \times N(0,0,0,0,2,3)(y) \}_{i=1}^5,
    \]
    with $(x^1_1, \dots, x^1_{10}) = (4,6,8,10,12,14,14,14,14,14)$, and
    \[
    \{\mathcal{B}^{\mathrm{ray}}_i\}_{i=6}^{10} = \{ N(x^2_i, \dots, x^2_{i+5})(x) \times N(7,8,10,10,10,10)(y) \}_{i=1}^5,
    \]
    with $(x^2_1, \dots, x^2_{10}) = (0,0,0,0,0,2,4,6,8,10)$.
    \item \textbf{Tensor product mesh}: The remaining 99 local tensor product B-splines, denoted as $\{\mathcal{B}^{\mathrm{cross}}_i\}_{i=1}^{99}$, form the tensor product component.
\end{itemize}
These 110 PT-spline basis functions form a basis for $S_{4,4}(\mathscr{T})$.
\end{example}

\subsection{\textcolor{blue}{Surface approximation and isogeometric analysis with PT-splines}}
\textcolor{blue}{In this section, we apply PT-splines constructed above for surface approximation and numerical solution of partial differential equations.}

\begin{figure}[htbp]  
    \centering    %
    \subfigure[] %
    {
        \begin{minipage}[t]{0.2\textwidth}
            \centering      %
            \includegraphics[width=1\textwidth]{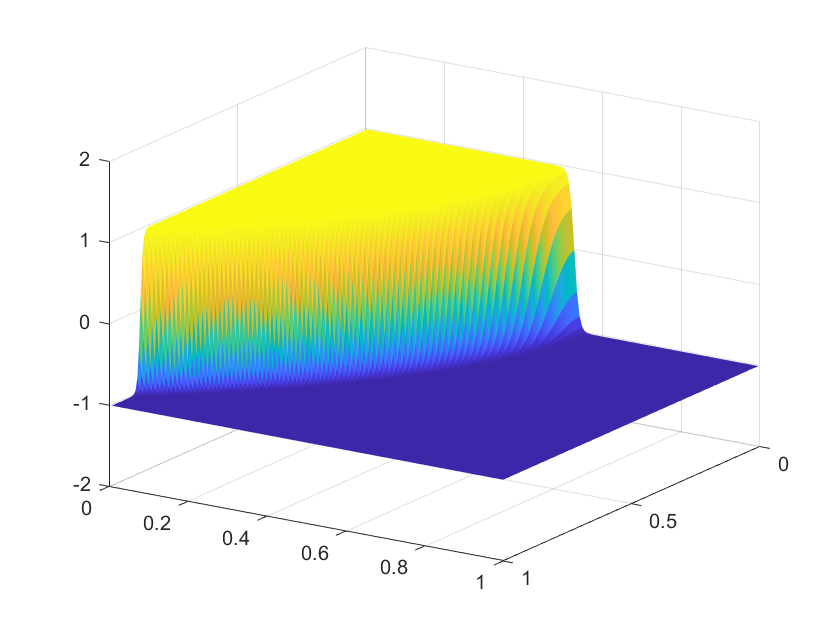}

        \end{minipage}
    }%
    \subfigure[] %
    {
        \begin{minipage}[t]{0.4\textwidth}
            \centering      %
            \includegraphics[width=0.48\textwidth]{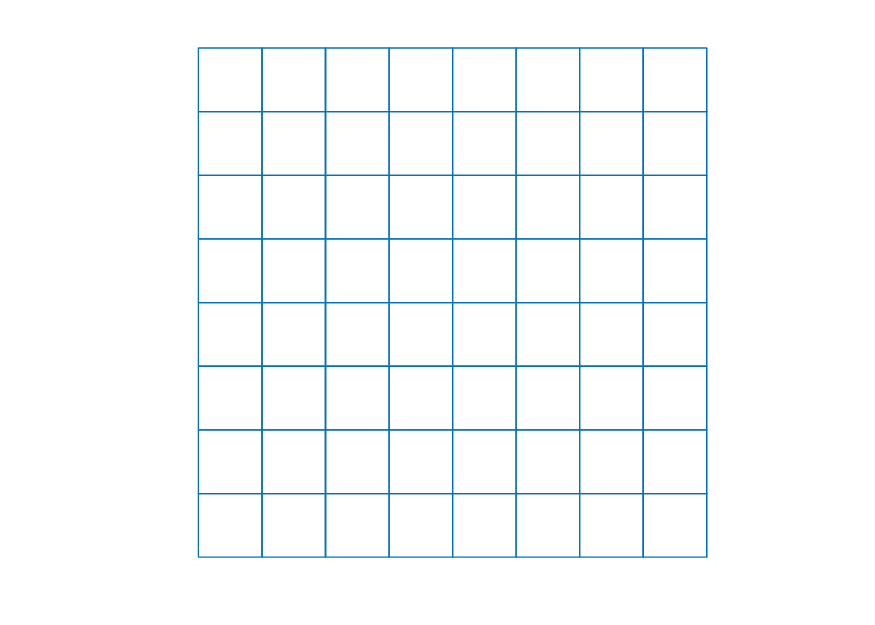} 
            \includegraphics[width=0.48\textwidth]{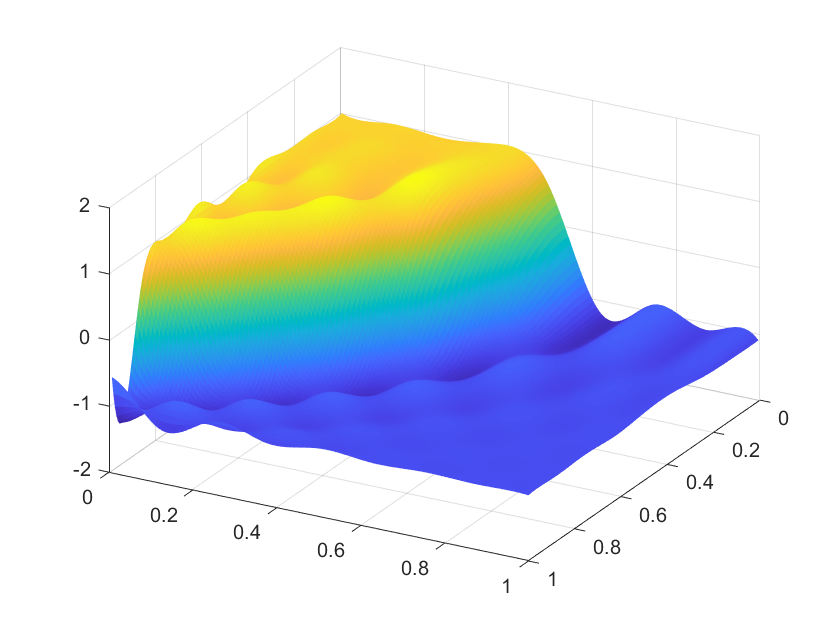}
            \end{minipage}
    }%
    \subfigure[] %
    {
        \begin{minipage}[t]{0.4\textwidth}
            \centering      %
            \includegraphics[width=0.48\textwidth]{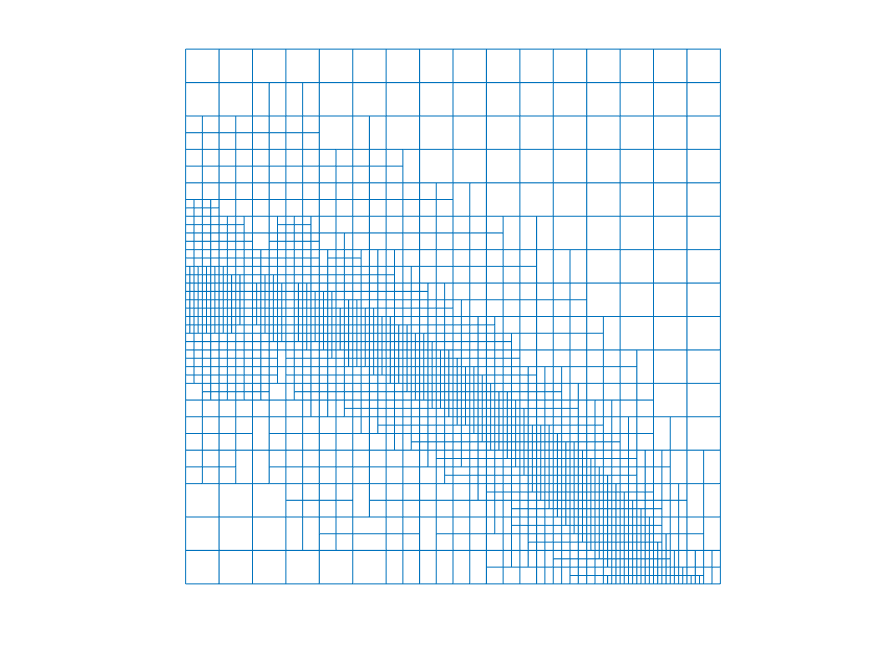} 
            \includegraphics[width=0.48\textwidth]{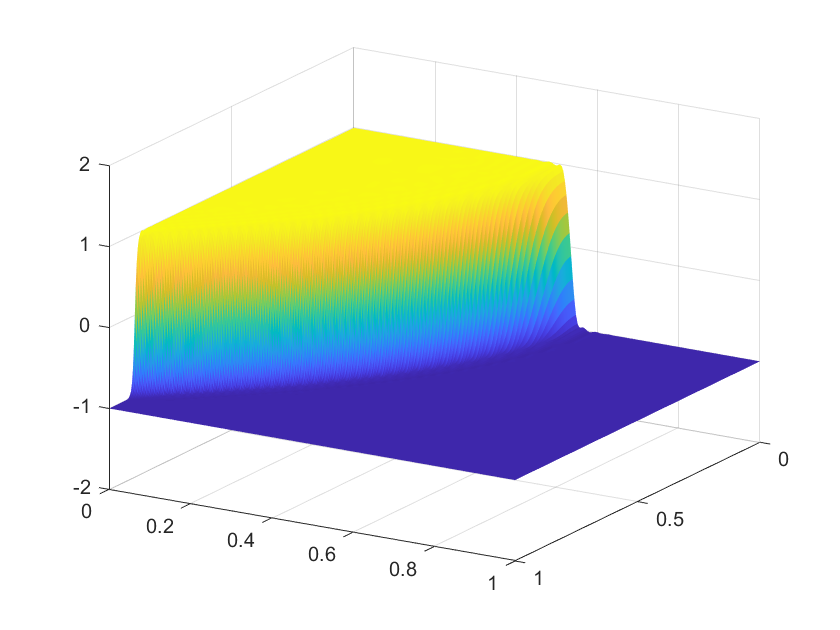}
            \end{minipage}
    }%
    
    \captionsetup{font={scriptsize}}
    \caption{The exact solution and fitting results using PT-spline basis in Example~\ref{exm sur app}.} 
    \label{sfex}  %
\end{figure}

\begin{example}\label{exm sur app}
\textcolor{blue}{Figure~\ref{sfex}(a) illustrate the surface of $$u(x,y)=-\mathrm{tanh}\frac{0.9+\sqrt{x^2+1.5y+0.01}}{0.01}.$$ 
}
\textcolor{blue}{We employ the least squares method, selecting \(151 \times 151\) uniformly distributed points in \([0, 1] \times [0, 1]\) as fitting points to approximate \(u\). The initial grid is an \(8 \times 8\) tensor product mesh, shown in Figure~\ref{sfex}(b) with its corresponding fitting surface. After refinement by adding horizontal or vertical lines, we obtain the T-mesh and fitting functions depicted in Figure~\ref{sfex}(c).}
\end{example}

\begin{figure}[htbp]  
    \centering    %
    \subfigure[] %
    {
        \begin{minipage}[t]{0.3\textwidth}
            \centering      %
            \includegraphics[width=1\textwidth]{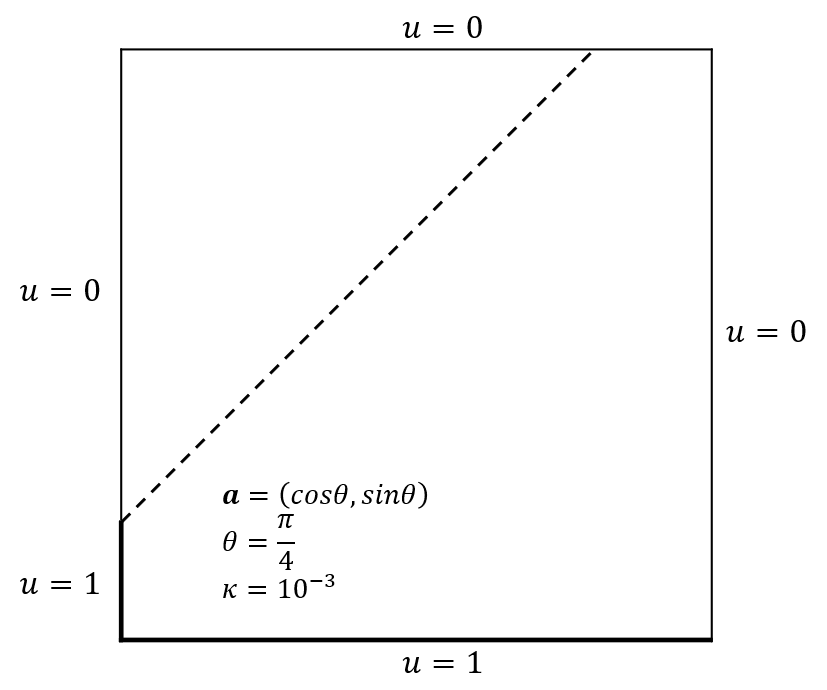}

        \end{minipage}
    }%
    \subfigure[] %
    {
        \begin{minipage}[t]{0.7\textwidth}
            \centering      %
            \includegraphics[width=0.48\textwidth]{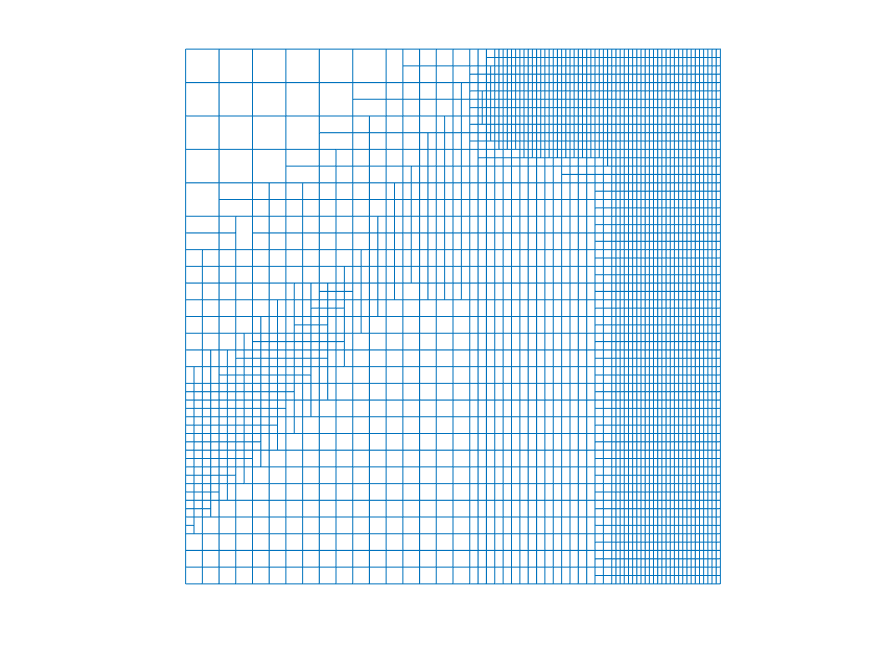} 
            \includegraphics[width=0.48\textwidth]{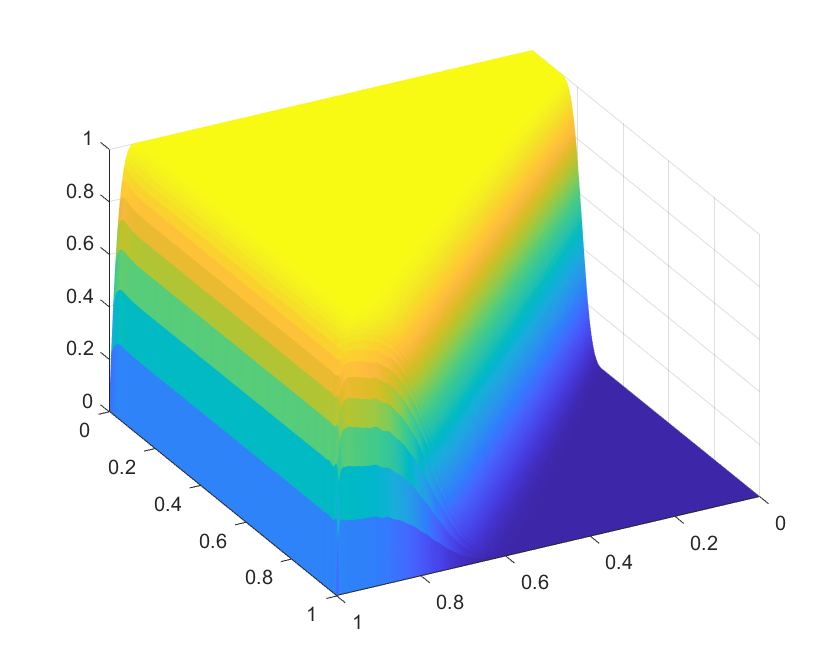}
            \end{minipage}
    }%
    
    \captionsetup{font={scriptsize}}
    \caption{The numerical solution using PT-spline basis in Example~\ref{exm phat pde}.} 
    \label{pdeex}  %
\end{figure}

\begin{example}\label{exm phat pde}
\textcolor{blue}{We use PT-spline basis to compute the numerical solution of linear advection-diffusion equation\begin{equation}
\begin{aligned}
 \mathbf{a}\cdot \nabla u-\nabla\cdot(\mathbf{\kappa}\nabla u) &=f\ \ \ in\ \ \ \Omega,\\
 u&=g\ \ \ on\ \ \ \Gamma_D   
\end{aligned}
\end{equation} where $\Omega $ is the unit suqare $[0,1]\times [0,1]$ and the relevant parameters in our example are shown in figure~\ref{pdeex}(a). We use streamline upwind/Petrov-Galerken stabilization and the final solution and the T-mesh are shown in figure~\ref{pdeex}(b).}
\end{example}

\section{Conclusion}
This paper presents a method for constructing basis functions for polynomial spline spaces with maximal smoothness over arbitrary T-meshes $S_{d_1,d_2}(\mathscr{T})$. This is a novel method to construct basis functions for spline spaces over arbitrary T-meshes. The proposed PT-spline basis construction involves two primary steps: (1) defining basis functions for the spline space over an extended diagonalizable T-mesh, based on the dimension formula for such meshes (Eq.~\eqref{dim diag simplify}); and (2) applying Extended Edge Elimination (EEE) conditions to connect the basis functions of the extended T-mesh to those of the original T-mesh. The PT-spline basis ensures linear independence and completeness. An algorithm for its construction is provided and compared with LR B-splines, a known spline formulation for T-meshes \cite{dokken2013polynomial}. Unlike PT-splines, LR B-splines lack linear independence, limiting their suitability as a basis for T-mesh spline spaces, and are restricted to LR-meshes, a specific T-mesh subtype. In contrast, PT-splines are applicable to any T-mesh. 

Future work may explore reducing PT-spline basis function support, optimizing algorithms to minimize added edges, ensure partition of unity, and applying PT-splines in isogeometric analysis and other domains. The method extends to polynomial spline spaces with arbitrary order of smoothness is worthy of future research.

\section{Acknowledgement}
We thank Professor Xin Li from the University of Science and Technology of China for providing an example of a diagonalizable T-mesh with T-cycles (Example~\ref{T-cycle}), significantly enhancing the applicability of this study. We also express gratitude to the anonymous reviewers for their valuable comments and suggestions. This work is supported by NSFC of China (no. 12494550 and 12494555).

\bibliographystyle{IEEEtran}

\bibliography{template} 


\end{document}